\documentclass[11pt]{article}
\usepackage{amsfonts}
\usepackage{graphics,latexsym,amssymb,amsmath,amscd}
\RequirePackage{graphics}
\usepackage{graphicx}
\usepackage{psfrag}
\input epsf
\def\Im{\hbox{\rm Im}\,}
\def\[#1\]{\begin{eqnarray*}#1\end{eqnarray*}}
\def\Re{\hbox{\rm Re}\,}

\def\Heis{{\mathfrak N}}

\def\phi{\varphi}

\newtheorem{thm}{Theorem}[section]
\newtheorem{dfn}[thm]{Definition}
\newtheorem{rmk}[thm]{Remark}
\newtheorem{cor}[thm]{Corollary}
\newtheorem{prop}[thm]{Proposition}
\newtheorem{lem}[thm]{Lemma}

\newcommand{\Pf}{{\sc Proof}. }
\newcommand{\EPf}{\hbox{}\hfill$\Box$\vspace{.5cm}}

\input xy
\xyoption{all}

\newcommand{\C}{{{\mathbb C}}}
\newcommand{\R}{{{\mathbb R}}}
\newcommand{\Z}{{{\mathbb Z}}}
\newcommand{\Q}{{{\mathbb Q}}}
\newcommand{\PP}{{{\mathbb P}}}
\date{}
\begin{document}
\title{A combinatorial invariant for Spherical CR structures}
\author{Elisha Falbel\\
        Institut de Math\'{e}matiques\\
        Universit\'{e} Pierre et Marie Curie\\
        4, Place Jussieu\\
        F-75252, Paris, France\\
        E-mail:falbel@math.jussieu.fr
        \and
        Qingxue Wang \thanks{Q.Wang was supported by NSFC grant \#10801034.}\\
        School of Mathematical Sciences \\
        Fudan  University\\
        Shanghai, 200433\\
        P.R. China\\
        E-mail:qxwang@fudan.edu.cn}
\maketitle

\begin{abstract}
We study a cross-ratio of four generic points of $S^3$
which comes from spherical CR geometry. We construct a
homomorphism from a certain group generated by generic
configurations of four points in $S^3$
to the pre-Bloch group $\mathcal {P}(\C)$. If $M$ is
a $3$-dimensional spherical CR manifold with a CR triangulation, by
our homomorphism, we get a $\mathcal {P}(\C)$-valued invariant for
$M$. We show that when applying to it the Bloch-Wigner function, it is zero.
Under some conditions on $M$, we show the invariant lies in the Bloch group $\mathcal B(k)$, where $k$ is the field
generated by the cross-ratio. For a CR triangulation of Whitehead link complement, we show its invariant is a
non-trivial torsion in $\mathcal B(k)$.
\end{abstract}

\section{Introduction}

One can define a number of invariants out of a
 3-manifold $M$ equipped with a complete real hyperbolic structure.
For instance,
 by Mostow's rigidity theorem, the volume of a hyperbolic 3-manifold
turns out to be a topological invariant
of the underlying manifold.  Another such invariant is Chern-Simons invariant $CS(M)$
with values in $\R/\Z$ (\cite{CS}).
Both invariants can be seen as arising from an invariant
associated to a hyperbolic manifold with values in the Bloch group (for the definition see Section
\ref{section:preliminaries}). The Bloch group $\mathcal {B}(\C)$
is a subgroup of the pre-Bloch group $\mathcal {P}(\C)$ which is defined
as the abelian group generated by all the points in $\C\setminus\{0,1\}$ quotiented by the 5-term relations.  The
volume and the Chern-Simons invariant
can then be seen through a function (the Bloch regulator)
$$
 \mathcal {B}(\C)\rightarrow {\C}/{\Q}.
$$
The imaginary part being the volume and the real part being $CS(M)$ mod $\Q$.

Another geometric structure on 3-manifolds
which has been studied for a long time is the Cauchy-Riemann (CR) structure.
  More precisely,
consider $S^3\subset \C^2$ with the contact structure obtained as the intersection $D=TS^3\cap JTS^3$ where $J$ is the
multiplication by $i$ in $\C^2$.
The operator $J$ restricted to $D$ defines the standard CR structure on $S^3$.
The group of CR-automorphisms of $S^3$ is $PU(2,1)$ and we say
that a manifold $M$ has a spherical CR structure if it has a $(S^3, PU(2,1))$-
geometric structure.

A configuration of four points in $S^3$ can be thought
as defining a CR simplex which can be parametrised through certain cross-ratios.
In this  paper we obtain an invariant associated to a simplicial complex
by CR simplices
with values in the pre-Bloch group ${\mathcal P}(\C)$ and
with some additional hypothesis, it lies
in $\mathcal B(\C)$.  In fact one can define the invariant
in ${\mathcal P}(k)$ where $k$ is the field generated by the cross-ratios
of the simplicial complex.  In the case $k$ is an imaginary quadratic extension
of a totally real field and if the invariant is in $\mathcal B(k)$, a theorem by Borel will imply that the element is
torsion and therefore $CS(M)=0$.
It is interesting to compare this result with the real
hyperbolic geometry case; Neumann and Yang \cite{NY1}
proved with the same hypothesis that
$CS(M)=0$, although the corresponding element in the Bloch group is never a torsion since the hyperbolic volume is
non-vanishing. On the other hand, in Section \ref{whitehead}, we have  examples with  invariants which are non-trivial
torsions in the Bloch group.
%It is interesting to have quite different behaviors of the invariants in the Bloch group
%for real hyperbolic geometry and spherical CR geometry.

The paper is organized as follows.
In the second section, we recall the basic definitions of (pre)-Bloch group, cross-ratio structures, triangulations,
and CR geometry.
In the third section, we study a complex of configurations of generic
 points in $S^3$, and a homomorphism to the pre-Bloch group $\mathcal {P}(\C)$.
 For a $3$-dimensional spherical CR manifold $M$ with a given CR triangulation, we associate  a $\mathcal
 {P}(\C)$-valued invariant $[M]$
which is independent of the triangulation. We show that when applying
to it  the Bloch-Wigner dilogarithm function, it is zero. In the fourth
section, we show that when  $M$ is closed with holonomy with
coefficients in a number field or non-compact with unipotent
parabolic boundary holonomy,
 the invariant $[M]$ lies in the Bloch group $\mathcal B(k)$, where $k$ is the field
generated by the cross-ratio.
For such an $M$, we define its Chern-Simons invariant as the real part
of $\rho([M])$, where $\rho$ is the Bloch regulator map.
In the last section, we compute the invariant for certain simplicial CR structures on the complement of the figure
eight knot and the Whitehead link.  The torsion element we find associated to the Whitehead link and figure-eight knot (see Section 5) is a very basic torsion element and it is null if we adopt a different definition of the Bloch group (see Remark after Definition \ref{invariant}).  On the other hand, we do not have examples of CR simplicial complexes
with non-torsion invariant in the Bloch group.

\section{Preliminaries}\label{section:preliminaries}
In this section, we recall the basic definitions and
 properties of (pre)-Bloch group, cross-ratio structures, triangulations,
and CR geometry.

\subsection{The Bloch group}

We consider an arbitrary field $F$ in the next definition
although we will only use $\C$, $\R$ and  number fields
in this paper. There are several definitions of the Bloch group in the literature but we
 will follow conventions of
\cite{S1}.
\begin{dfn}
The pre-Bloch group $\mathcal{P}(F)$ is the quotient of the free
abelian group $\Z[F\setminus\{0,1\}]$ by the subgroup generated by
the 5-term relations
\begin{equation}\label{5term}
 [x]-[y]+[\frac{y}{x}]-[\frac{1-x^{-1}}{1-y^{-1}}]+[
  \frac{1-x}{1-y}],\ \  \forall\, x,y\in F\setminus\{0,1\}.
\end{equation}
\end{dfn}

Consider the tensor product $F^*\otimes_{_\Z} F^*$, where $F^*$ is the
multiplicative group of $F$. It is an abelian group
satisfying, for $n$ an integer, $n(a\otimes b)= a^n\otimes b= a\otimes b^n$.

Let $T=\langle x\otimes y+y\otimes x\ |\ x,y\in F^{*} \rangle$ be the subgroup of $F^{*}\otimes_{_\Z}F^{*}$ generated
by $x\otimes y+y\otimes x$, where $x,y\in F^{*}$.
\begin{dfn}
$\bigwedge^2F^{*}= (F^{*}\otimes_{_\Z}F^{*})/T$. For $x,y\in F^{*}$, we will denote by $x\wedge y$ the image of
$x\otimes y$ in $\bigwedge^2F^{*}$.
\end{dfn}
Note that, for $x,y \in F^{*}$, we have $x\wedge y=-y\wedge x$ and $2\, x\wedge x=0$.
But $x\wedge x=0$ is
not necessarily true in $\bigwedge^2F^{*}$.
\begin{dfn}
The Bloch group $\mathcal{B}(F)$ is the kernel of the homomorphism
$
\delta : \mathcal{P}(F) \rightarrow \bigwedge^{2} F^{*},
$
which is defined on  generators of $\mathcal{P}(F)$  by
$\delta ([z])=z\wedge(1-z)$.
\end{dfn}

When $F=\C$, $\mathcal{P}(\C)$ (\cite{DS} theorem 4.16) and  $\mathcal{B}(\C)$ (\cite{S1})
are uniquely divisible groups and, in fact, are $\Q$-vector
spaces with infinite dimension (\cite{S1}).
In particular they have no torsion.  On the other hand, when the field is $\R$
there exists torsion.  In particular, for all $x\in \R-\{0,1\}$,
the element $[x]+[1-x]\in \mathcal{B}(\R)$ does not depend on $x$ and
has order six (see \cite{S1} prop. 1.1 pg 220).

We will need the following result for $\mathcal{B}(k)$ when $k$ contains
a cubic root of unity:

\begin{lem}\label{2beta}
Suppose $k$ is a field containing a cubic root of unity.  Then
$$
2[z]=2\left [ \frac{1}{1-z}\right ], \; \forall z\in k\setminus\{0,1\}.
$$
\end{lem}
\Pf By \cite{S1} Lemma 1.5 (c), we have
$$
  2([z]+[1-z])=0.
$$
From Lemma 1.2 of \cite{S1}, we obtain
$$
2([1-z]+\left [ \frac{1}{1-z}\right ])=0.
$$
Therefore,
$$
2[z]=-2[1-z]=2\left [\frac{1}{1-z}\right ].
$$
\EPf

Consider the complex conjugation in $\C$ and its extension
to an involution:
$$
\sigma :\Z[\C\setminus \{0,1\}]
\rightarrow \Z[\C\setminus \{0,1\}].
$$  As $\sigma$
preserves the 5-term relations, it induces an involution on the
pre-Bloch group $\mathcal{P}(\C)$ which we will also denote by $\sigma$. Let
$$
 \mathcal{P}(\C)=\mathcal{P}(\C)^{+}+\mathcal{P}(\C)^{-}
$$
be the decomposition of $\mathcal{P}(\C)$ into the two subgroups preserved by
the involution. They are the eigenspaces of $\sigma$ acting on the
$\Q$-vector space $\mathcal{P}(\C)$. We have $\sigma(z)=z$ for $z\in \mathcal{P}(\C)^{+}$ and
 $\sigma(z)=-z$ for $z\in \mathcal{P}(\C)^{-}$.  Analogously
$$
 \mathcal{B}(\C)=\mathcal{B}(\C)^{+}+\mathcal{B}(\C)^{-}.
$$

By \cite[~Lemma 1.3]{S1}, for any $x,y\in F-\{0,1\}$, $[x]+[1-x]=[y]+[1-y]$ in $\mathcal {P}(F)$. That is, the element
$[x]+[1-x]\in \mathcal{B}(F)$ is independent of the choice of $x$, hence the following definition makes sense.

\begin{dfn}\label{bcf}
$c_F:=[x]+[1-x]\in \mathcal{B}(F)$, where $x$ is any given element of $F \setminus \{0,1\}$.
\end{dfn}

Let $K^{M}_{*}(F)$ denote the Milnor K-groups of $F$, and $K_{*}(F)$ denote the Quillen's algebraic K-groups of $F$.
It is well-known that $K^{M}_1(F)\cong K_{1}(F)=F^{*}$ and $K^{M}_2(F)\cong K_{2}(F)$. There is a natural map from
$K^{M}_{n}(F)$ to $K_{n}(F)$ for each $n$. The cokernel of this map is called the group of indecomposable elements,
denoted by $K_{n}^{\text{ind}}(F)$.

For $n=3$ and $F$ an infinite field, we have the following fundamental exact sequence, due to Suslin (\cite[~Theorem
5.2]{S1}):
\begin{equation}\label{sus1}
   \begin{CD}
   0 @>>>\text{Tor}(\mu(F),\mu(F))^{\sim} @>>> K_{3}^{\text{ind}}(F) @>s_F>> \mathcal{B}(F) @>>> 0, \\
   \end{CD}
\end{equation}
where $\mu(F)$ is the group of roots of unity in $F$. If $\text{char}F=2$,
$\text{Tor}(\mu(F),\mu(F))^{\sim}=\text{Tor}(\mu(F),\mu(F))$. If $\text{char}F\ne 2$,
$\text{Tor}(\mu(F),\mu(F))^{\sim}$ is the unique nontrivial extension of $ \text{Tor}(\mu(F),\mu(F))$ by $\Z/2$, that
is, we have the nontrivial extension
$$
 \begin{CD}
   0 @>>> \Z/2@>>> \text{Tor}(\mu(F),\mu(F))^{\sim} @>>> \text{Tor}(\mu(F),\mu(F)) @>>> 0, \\
   \end{CD}
$$

For a field extension $E/F$, there are natural homomorphisms $\alpha:K_{3}^{\text{ind}}(F)\rightarrow
K_{3}^{\text{ind}}(E)$ and $\beta: \mathcal{B}(F)\rightarrow \mathcal{B}(E)$.
\begin{lem}\label{blem1}
Let $F$ be an infinite field and $E/F$  a field extension. The homomorphism \\$s_F:K_{3}^{\text{ind}}(F)\rightarrow
\mathcal{B}(F)$ in (\ref{sus1}) is functorial in $E/F$ with respect to $\alpha, \beta$. That is, we have the following
commutative diagram:
$$
 \begin{CD}
  K_{3}^{\text{ind}}(E)@>s_E>>\mathcal{B}(E)\\
  @A\alpha AA @A\beta AA\\
  K_{3}^{\text{ind}}(F)@>s_F>>\mathcal{B}(F)
 \end{CD}
$$
\end{lem}
\Pf
It follows from the construction of the exact sequence (\ref{sus1}) in \cite{S1}. The first step in the construction is the homomorphism
from $H_3(GL(F))$ to $\mathcal{B}(F)$ via spectral sequences, see \cite[~Proposition 3.1, Theorem 4.1]{S1}. The second
step is the Hurewicz homomophism from $K_3(F)$ to $H_3(GL(F))$, see \cite[~Lemma 5.4, Theorem 5.1]{S1}. It is clear
that both steps are functorial with respect to the field extension $E/F$. Hence the diagram commutes.
\EPf

\begin{lem}\label{blem2}
Let $E$ be a field of characteristic $0$. If $\text{Tor}(\mu(E),\mu(E))$ has no element of order $3$, then $2c_E\in
B(E)$ has order $3$.
\end{lem}
\Pf
Denote by  $o(g)$ the order of a group element $g$. Consider the field extension $E/\Q$. By
\cite[~Corollary 5.3]{S1}, the Bloch group
 $\mathcal{B}(\Q)$ is generated by $c_{\Q}$ and $o(c_{\Q})=6$. Hence $o(2c_{\Q})=3$. By Lemma \ref{blem1}, we have the
 commutative diagram:
$$
 \begin{CD}
  K_{3}^{\text{ind}}(E)@>s_E>>\mathcal{B}(E)\\
  @A\alpha AA @A\beta AA\\
  K_{3}^{\text{ind}}(\Q)@>s_{\Q}>>\mathcal{B}(\Q)
 \end{CD}
$$
By definition $\beta(2c_{\Q})=2c_E$. Since $o(2c_{\Q})=3$, $o(2c_E)=1$ or $3$. It suffices to show that $2c_E\ne 0$.
Since $\mu(\Q)=\Z/2$, $\text{Tor}(\Z/2,\Z/2)=\Z/2$. Hence $\text{Tor}(\mu(\Q),\mu(\Q))^{\sim}=\Z/4$. Now by the
fundamental exact sequence ($\ref{sus1}$), there is a nonzero element $q\in K_{3}^{\text{ind}}(\Q)$ such that
$s_{\Q}(q)=2c_{\Q}$. Since $o(2c_{\Q})=3$, we see that $3|o(q)$. Let $r=\alpha(q)$. By the above commutative diagram,
we get $s_E(r)=2c_E$. Therefore, it suffices to show that $r\notin \text{ker}\,s_{E}$. Indeed, since
$\text{Tor}(\mu(E),\mu(E))$
has no element of order $3$, we see that $\text{Tor}(\mu(F),\mu(F))^{\sim}$
does not have an element of order $3$ either. Hence it has no element of order
divided by $3$. By \cite[~Proposition 11.3]{S2}, the homomorphism $\alpha$ is injective. Therefore, $3|o(r)=o(q)$. Now
the fundamental exact sequence ($\ref{sus1}$) implies that $r\notin \text{ker}\,s_{E}$. Hence $o(2c_E)=3$.
\EPf

Consider now a number field $k$, that is, an extension of $\Q$ of
degree $d<\infty$. We
have that $k=\Q(t)$, where $t\in \C$ satisfies an irreducible
polynomial  with coefficients in $\Q$ with degree $d$.
Each root of the irreducible
polynomial determines a field embedding  $\sigma : k\rightarrow \C$.
There are
$r_1$ real embeddings (when $\sigma(k)\subset \R$) and $r_2$
pairs of complex conjugate embeddings with $d=r_1+2r_2$.

\begin{dfn}\label{BWf} The Bloch-Wigner function is
$$
  D(x)=\arg{(1-x)}\log{|x|}-\Im (\int_{0}^{x}\log{(1-t)}\frac{dt}{t}),
$$
\end{dfn}
It is well-defined and real analytic on $\C-\{0,1\}$ and
extends to a continuous function on $\C P^1$ by
defining $D(0) = D(1) = D(\infty) = 0$. It is well-known
that it satisfies the 5-term relation. Hence it gives rise
to a well-defined map:
$$
  D: \mathcal{P}(k)\rightarrow \R,
$$
given by
$$
 D(\sum_{i=1}^{k}n_i[x_i])=\sum_{i=1}^{k}n_iD(x_i).
$$
Generalizing Dirichlet's units theorem, Borel \cite{Bo} proved the
following description
of $\mathcal{B}(k)$.  His results are more general and describe
the higher K-theory of number fields.  The relation to the Bloch group and
the use of the dilogarithm is due to work of Dupont, Sah, Bloch and Suslin.
In the following we let $\sigma_i, \bar \sigma_i$, $1\leq i\leq r_2$,
be the $r_2$ pairs of complex embeddings.

\begin{thm}[\cite{Bo}] Consider the map
$r :  \mathcal{B}(k)\rightarrow
\R^{r_2}$ given by
$$
[z]\rightarrow \left ( D(\sigma_1(z)),\cdots ,D(\sigma_{r_2}(z))\right ).
$$
Then $Im (r)$ is a lattice in $\R^{r_2}$ and $Ker(r)$ is the torsion subgroup
of $\mathcal{B}(k)$.
\end{thm}

Recall that we have the Bloch regulator map:
$$
\rho : \mathcal{B}(\C)\rightarrow \C/\pi^2\Q=\R/\pi^2\Q \oplus i \R.
$$
It is known that the imaginary part of $\rho$ coincides with the
Bloch-Wigner function $D$.

\subsection{Tetrahedra with cross-ratio structures}
In the following   we recall the definition of cross-ratio
structures of \cite{F3}.

Consider a set of four elements $\Delta= \{ p_0,p_1,p_2,p_3 \}$.  We
call  $p_i$, $0\leq i\leq 3$ the vertices of $\Delta$.  Let
$O\Delta$ be the set of all orderings of $\Delta$.  We will denote
an element of $O\Delta$ by $[p_i,p_j,p_k,p_l]$ (where $\{
i,j,k,l\}=\{0,1,2,3\}$) and call it a simplex although we only deal
with configurations of four points. Given $\Delta$, there are 24
simplices divided in two classes $O\Delta^+$ (containing
$[p_0,p_1,p_2,p_3]$)
 and $O\Delta^-$ (containing $[p_0,p_1,p_3,p_2]$)
  of 12 elements each.  Each class is an orbit of the even permutation group acting on $O\Delta$.

The following definition assigns similarity invariants
to each vertex of a configuration of four points.

\begin{dfn} A cross-ratio structure on a set of four points  $\Delta=\{ p_0,p_1,p_2,p_3 \}$ is a function defined on
the ordered quadruples
$$
{\bf X} : O\Delta\rightarrow \C\setminus \{0,1\}
$$
satisfying, if
$(i,j,k,l)$ is any permutation of $(0,1,2,3)$, the relations
\begin{enumerate}
\item
$$
{\bf X}(p_i,p_j,p_k,p_l)=\frac{1}{{\bf X}(p_i,p_j,p_l,p_k)}.
$$
\item  (similarity relations)
$$
{\bf X}(p_i,p_j,p_k,p_l)=\frac{1}{1-{\bf X}(p_i,p_l,p_j,p_k)}.
$$
\end{enumerate}
\end{dfn}

{\bf Remarks}

\begin{figure}
\setlength{\unitlength}{1cm}
\begin{center}
\begin{picture}(7,10)
\psfrag{z1}{$z_{01}$}
\psfrag{z1'}{$z_{10}$}
\psfrag{z~1}{$z_{23}$}
\psfrag{z~1'}{$z_{32}$}
\psfrag{z2}{$z_{02}$}
\psfrag{z2'}{$z_{13}$}
\psfrag{z~2}{$ z_{20}$}
\psfrag{z~2'}{$ z_{31}$}
\psfrag{z3}{$z_{03}$}
\psfrag{z3'}{$z_{12}$}
\psfrag{z~3}{$ z_{21}$}
\psfrag{z~3'}{$ z_{30}$}
\psfrag{q1}{$p_2$}
\psfrag{q2}{$p_3$}
\psfrag{p1}{$p_0$}
\psfrag{p2}{$p_1$}
{\scalebox{.8}{\includegraphics[height=8cm,width=8cm]{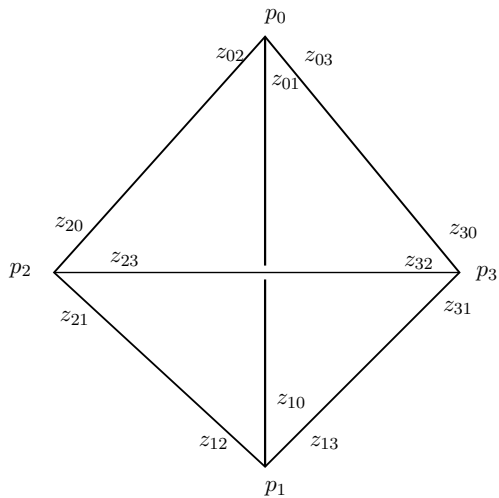}}}
\end{picture}
\end{center}
\caption{\sl Parameters for a cross-ratio structure}\label{Figure:ptetrahedron}
\end{figure}

\begin{enumerate}
\item To visualize the definition we refer to Figure
 \ref{Figure:ptetrahedron}.
For each $[p_i,p_j,p_k,p_l]\in O\Delta^+$ we define
$$
z_{ij}={\bf X}(p_i,p_j,p_k,p_l).
$$
We interpret $z_{ij}$ as a cross-ratio associated to the edge $[ij]$ at the vertex $i$.
Cross-ratios of elements of  $O\Delta^-$ are obtained taking inverses by the first
symmetry.
In the following we shall denote by a sequence of numbers $(ijkl)$ the
corresponding invariant ${\bf X}(u_i,u_j,u_k,u_l)$.
\item The similarity relations can be used to reduce the number of variables to four,
one for each vertex.  One can use, for instance,
$(z_{01},z_{10},z_{23},z_{32})\in
\left ( \C\setminus \{0,1\}\right )^{4}$.
\item
If we impose that
$$
z_{12}=z_{21}=z_{34}=z_{43},
$$
we may interpret the configuration of four points as a configuration
of points in $\C P^1$ with $z_{12}$ one of the cross-ratios.
\end{enumerate}

\subsection{Triangulations, cross-ratio structures and an invariant}

A triangulation of a manifold is an explicit description of it as a
simplicial complex. Equivalence of triangulations was studied in \cite{A}.  Alexander defined
certain moves on triangulations and showed that any two triangulations of a
closed manifold are related by a sequence of these moves.
A particularly simple description of that equivalence was obtained by
Pachner \cite{P} through some elementary moves.  In particular, Pachner
proved that any two triangulations are obtained from each other
 through the following two
moves and their inverses.

 \begin{dfn} An elementary (or Pachner) move
 in a simplicial complex defined by a triangulation
is
\begin{itemize}
\item (2-3 move) the substitution of two simplices
 $[u_0,u_1,u_2,u_3]-[u_0,u_1,u_2,u_4]$, which have a common face, by
three simplices $-[u_0,u_1,u_3,u_4]+[u_0,u_2,u_3,u_4]-[u_1,u_2,u_3,u_4]$
with one common edge and vice-versa.
\item (1-4 move) the substitution
of a simplex $[u_0,u_1,u_2,u_3]$ by four simplices
$$
[u_0,u_1,u_2,u_4]-[u_0,u_1,u_3,u_4]+[u_0,u_2,u_3,u_4]-[u_1,u_2,u_3,u_4]
$$ by adding a new vertex $u_4$ and vice-versa.
\end{itemize}
\end{dfn}

To deal with manifolds with boundary it is useful to
introduce ideal triangulations.  It is important in the following to consider
singular triangulations, meaning that we allow self or multiple intersections
between simplices (along 2-faces).

\begin{dfn}
Let $M$ be a compact manifold with boundary $\partial M$.  Denote by $\hat M$
the manifold obtained by collapsing each boundary component.
 An ideal
triangulation is a singular triangulation of $\hat M$ such that its vertices
are  the points of $\hat M$ obtained from collapsing the boundary components
of $M$.
\end{dfn}

A related concept we will not deal with in this paper is that of a spine
of a 3-manifold with boundary.
Ideal triangulations are in natural bijection to standard spines and existence
and equivalence of ideal triangulations were proven through the corresponding
existence and equivalence theorems for standard spines (a general reference is
\cite{Ma1}).
We will need the following theorem (\cite{Ma,P}).

\begin{thm}\label{tr}
Any two ideal triangulations of a manifold with boundary can be obtained
from one another through the 2-3 moves.
\end{thm}

\begin{dfn}
Let $T$ be a (ideal, if the 3-manifold is not closed)
 triangulation of a 3-manifold.
%By this we
% mean a simplicial
%complex whose underlying topological space is a manifold if
%the vertices are deleted.
Let ${\bf X}(p_i,p_j,p_k,p_l)$ be a cross-ratio structure
 defined on the simplices. The pair $(T,{\bf X})$ is called
a cross-ratio structure associated to a triangulation
if the following compatibility conditions are satisfied:
\begin{enumerate}
\item Edge compatibility: If $[p_i,p_j,p_{m_0},p_{m_1}],
[p_i,p_j,p_{m_1},p_{m_2}],\cdots ,
[p_i,p_j,p_{m_n},p_{m_0}]$ are simplices having the  edge $[p_i,p_j]$
in common then
$$
{\bf X}(p_i,p_j,p_{m_0},p_{m_1})\cdots {\bf X}(p_i,p_j,p_{m_n},p_{m_0}) =1
$$
\item Face compatibility:  If $[p_i,p_j,p_k,p_l]$ and $[p_{i'},p_j,p_k,p_l]$
are two simplices with a common face $[p_j,p_k,p_l]$ then
$$
{\bf X}(p_j,p_i,p_k,p_l){\bf X}(p_k,p_i,p_l,p_j){\bf X}(p_l,p_i,p_j,p_k)
=
{\bf X}(p_j,p_{i'},p_k,p_l){\bf X}(p_k,p_{i'},p_l,p_j){\bf X}(p_l,p_{i'},p_j,p_k)
$$
\end{enumerate}
\end{dfn}

%{\bf Remark}

\begin{dfn}
To each cross-ratio structure associated to a triangulation as above we
define
the element
$$
\beta(T,{\bf X})= \sum_{s} ([z^s_{01}]+[z^s_{10}]+[z^s_{12}]+[z^s_{21}])\in {\mathcal P(\C)}
$$
where $s$ indexes the simplices in the triangulation and $z^s_{ij}$ are
the cross-ratios of the simplex $s$.
\end{dfn}

Consider $(T',{\bf X'})$ and $(T,{\bf X})$, two cross-ratio
structures on two triangulations obtained from each other by an
elementary move. Given the five vertices $\{ i,j,k,l,m\}=\{
0,1,2,3,4\}$, where the move is
concentrated as above, we write $(ijkl)$ the cross-ratio
defined by these points (which might be ${{\bf X'}(ijkl)}$ or ${{\bf
X}(ijkl)}$ according to which triangulation the simplex $[ijkl]$
belongs).

\begin{dfn}
We say that $(T',{\bf X'})$ is obtained from $(T,{\bf X})$ by an elementary
move if the triangulation $T'$ is obtained from $T$ by an elementary move and
  the following relations are satisfyed:
\begin{enumerate}
\item edge compatibility conditions
$$
(ijkl)=(ijkm)(ijml).
$$
\item face compatibility conditions
$$
(ijkl)(ljik)(kjli)=(imkl)(lmik)(kmli).
$$
\end{enumerate}
\end{dfn}

\begin{thm}\label{beta}  If $(T',{\bf X'})$ is obtained from $(T,{\bf X})$ by an elementary move
     then $\beta(T',{\bf X'})=\beta(T,{\bf X})$.
\end{thm}

\Pf It is exactly the proof of Theorem 5.2 in \cite{F3}.
\EPf

%{\bf Remarks}
%\begin{enumerate}
%\item
%if $(T,{\bf X})$ is the cross-ratio of a hyperbolic or CR triangulation we can
%obtain easily $(T',{\bf X'})$ from $(T,{\bf X})$ by an elementary move.  In fact
%a cross-ratio structure on a tetrahedron
% determines a configuration of four points in
%$\C P^1$ or $S^3$. Moreover, given five points in these spaces the cross-ratios
%of a combination of four points among them satisfy the compatibility conditions.
%\end{enumerate}

\subsection{CR geometry (see \cite{C,BS,G,J})}\label{section:CR}

\label{CRgeometry}

CR geometry is modeled on the {\sl Heisenberg group} $\Heis$,
 the set of pairs $(z, t)\in
{\C}\times{\R}$ with the product
$$
(z,t)\cdot (z',t') = (z+z', t + t' + 2 \Im z \overline{z}').
$$
The one point compactification of the Heisenberg group,
$\overline{\Heis}$, of $\Heis$ can be interpreted as $S^3$ which,
in turn, can be identified to the boundary of complex hyperbolic
space.

We consider the group $U(2,1)$ preserving the Hermitian form
$\langle z,w \rangle = w^*Jz$ defined on ${\mathbb C}^{3}$ by
the matrix
$$
J=\left ( \begin{array}{ccc}

                        0      &  0    &       1 \\

                        0       &  1    &       0\\

                        1       &  0    &       0

                \end{array} \right )
$$
and the following subspaces in ${\mathbb C}^{2,1}$ (that is, ${\mathbb C}^{3}$
with the above Hermitian form).
$$
        V_0 = \left\{ z\in {\mathbb C}^{3}-\{0\}\ \ :\ \
 \langle z,z\rangle = 0 \ \right\},
$$
$$
        V_-   = \left\{ z\in {\mathbb C}^{3}\ \ :\ \ \langle z,z\rangle < 0
\ \right\}.
$$
Let $\PP:{\C}^{3}\setminus\{ 0\} \rightarrow {\C}P^{2}$ be the
canonical projection.  Then
${\bf H}_{\C}^{2} = \PP(V_-)$ is the complex hyperbolic space and
$S^3=\partial{\bf H}_{\C}^{2} = \PP(V_0)$ can be identified to $\overline\Heis$.

The group of biholomorphic transformations of ${\bf H}_{\C}^{2}$ is then
$PU(2,1)$, the projectivization of $U(2,1)$.  It acts on $S^3$ by
CR transformations.
% An involution in $PU(2,1)$ has  a fixed point
%in the interior of complex hyperbolic space.  If it has fixed points in the
%boundary  of complex hyperbolic space, one shows that the set
%of fixed points is a topological circle, called $\C$-{\sl circle}.
  We
define $\C$-circles as boundaries of complex lines in ${\bf H}_{\C}^{2}$. Analogously,
$\R$-circles are boundaries of totally real totally geodesic two dimensional submanifolds in ${\bf H}_{\C}^{2}$.
Using the identification $S^3= \Heis \cup \{ \infty\}$,
one can define alternatively a $\C$-circle
 as any circle in $S^3$ which is obtained
from the vertical line $\{(0,t)\in \Heis \ |\ t\in \R\ \}\cup \{ \infty\}$ in the compactified Heisenberg
space by translation by an element of $PU(2,1)$.  Analogously, $\R$-circles are all obtained
by translations of the horizontal line $\{(x,0)\in \Heis \ |\ x\in \R\ \}\cup \{ \infty\}$.

A point $p=(z,t)$ in the Heisenberg group and the point $\infty$ are lifted
to the following points in $\C^{2,1}$:
$$
\hat{p}=\left[\begin{matrix}\frac{-|z|^2+it}{2} \\ z \\ 1 \end{matrix}\right]
\quad\hbox{ and }\quad
\hat{\infty}=\left[\begin{matrix} 1 \\ 0 \\ 0 \end{matrix}\right].
$$

Consider ${\mathbb C}^{2,1}$ with the Hermitian form defined by $J$.  The
Hermitian structure defines a Hermitian form on $\Lambda^3( {\mathbb C}^{2,1})$.
In fact,  let $\{ e_1, e_2, e_3\}$ be a basis of ${\mathbb C}^{2,1}$ such that
$$
\langle e_1, e_1\rangle=\langle e_2, e_2\rangle = -\langle e_3, e_3\rangle=1.
$$Then
 $e=e_1\wedge e_2\wedge e_3$ is a basis of $\Lambda^3( {\mathbb C}^{2,1})$ and
we define
$\langle e_1\wedge e_2\wedge e_3, e_1\wedge e_2\wedge e_3\rangle =-1$.

\begin{dfn}\label{boxproduct}
Let $v,w\in {\mathbb C}^{2,1}$.  We define the Hermitian cross-product
$v\boxtimes w$ by the formula
$$
\langle u, v\boxtimes w\rangle e = u\wedge v\wedge w,
$$
for all $u\in {\mathbb C}^{2,1}$.
\end{dfn}
In coordinates, $v=(v_0,v_1,v_2)$ and $w=(w_0,w_1,w_2)$, we compute
$$
v\boxtimes w=({\bar v}_0{\bar w}_1-{\bar v}_1{\bar w}_0,{\bar v}_2{\bar w}_0-{\bar v}_0{\bar w}_2,{\bar v}_1{\bar
w}_2-{\bar v}_2{\bar w}_1).
$$
So we observe that the coordinates of $v\boxtimes w$ are in the field
generated by the conjugates of the coordinates of $v$ and $w$.

\begin{dfn}
Given any three ordered points $p_0, p_1, p_2 \in \partial{\bf H}_{\C}^{2}$,
 we define Cartan's angular
invariant $\mathbb{A}$ as
$$
\mathbb{A}(p_0, p_1, p_2) =
\arg{(-\langle \hat{p}_0,\hat{p}_1\rangle \langle
\hat{p}_1,\hat{p}_2\rangle \langle \hat{p}_2,\hat{p}_0\rangle )}.
$$
\end{dfn}
Triple of points in $S^3$ are classified by their Cartan invariant according to the following proposition by Cartan
(see \cite{G}).
\begin{prop}\label{proposition:cartan}
Two ordered triples of pairwise distinct
 points  $(p_0, p_1, p_2)$ and $(p'_0, p'_1, p'_2)$
in $\partial{\bf H}_{\C}^{2}$ are equivalent under an element of $PU(2,1)$
if and only if
$$\mathbb{A}(p_0, p_1, p_2)=\mathbb{A}(p'_0, p'_1, p'_2).
$$
\end{prop}

\begin{dfn}\label{definition:generic}
We say that $n$ ($n\geq 1$)
points in $S^3$ are \emph{generic} if they are distinct
and, if $n\geq 3$,  any
three of them are not contained in a $\C$-circle.
\end{dfn}

As the action of $PU(2,1)$ is doubly transitive, a simple computation
gives the following normalization for a triple of generic points.
\begin{lem}
One can always normalize a triple
$(p_0,p_1,p_2)$ of generic points in $S^3$ so that in the
compactified Heisenberg group they are
$$
p_0=\infty\ \ \ p_1=(0,0)\ \ \ p_2=(1,t) \ \ \mbox{with} \ t\in \R .
$$
In that case $\tan \mathbb{A}(p_0, p_1, p_2)= t$.
\end{lem}

\subsection{Configurations of four points}

We refer to Figure \ref{Figure:ptetrahedron} to describe the
parameters of a tetrahedron (see also \cite{F2,F3}). Consider a
generic configuration of four (ordered) points $(p_0,p_1,p_2,p_3)$
in $S^3$. Fix one of them say $p_0$ and consider the projective
space of complex lines passing through it.  Then $p_1,p_2,p_3$
determine three points $t_1,t_2,t_3$ on $\C P^1$. The fourth point
corresponds to the complex line passing through $p_0$ and tangent to
$S^3$, call it $t_0$.  The cross-ratio of those four points in $\C
P^1$ is $z_{01}={\bf X}(t_0,t_1,t_2,t_3)$ (here, ${\bf X}$ is the
usual cross-ratio of four points in $\C P^1$).  We define
analogously the other invariants $z_{ij}$, $0\leq i,j \leq3$. If we
take $p_0=\infty$, the complex lines passing through $p_0$ intersect
$\Heis$ in vertical lines which are then determined by a coordinate
in $\C$. Up to Heisenberg translations, we can assume that
$p_1=(0,0)$ and $p_2=(1,s_2)$ and $p_3=(z_{01},s_3)$, $s_2,s_3\in
\R$. The corresponding points in $\C P^1$ will be
$\infty,0,1,z_{01}$. Therefore one ``sees'' at the vertex $p_0$ the
Euclidean triangle determined by $0,1,z_{01} \in \C$.

%\begin{figure}
%\setlength{\unitlength}{1cm}
%\begin{center}
%\begin{picture}(7,10)
%\psfrag{z1}{$z_1$}
%\psfrag{z1'}{$z_1'$}
%\psfrag{z~1}{$\tilde z_1$}
%\psfrag{z~1'}{$\tilde z_1'$}
%\psfrag{z2}{$z_2$}
%\psfrag{z2'}{$z_2'$}
%\psfrag{z~2}{$\tilde z_2$}
%\psfrag{z~2'}{$\tilde z_2'$}
%\psfrag{z3}{$z_3$}
%\psfrag{z3'}{$z_3'$}
%\psfrag{z~3}{$\tilde z_3$}
%\psfrag{z~3'}{$\tilde z_3'$}
%\psfrag{q1}{$p_3$}
%\psfrag{q2}{$p_4$}
%\psfrag{p1}{$p_1$}
%\psfrag{p2}{$p_2$}
%{\scalebox{.8}{\includegraphics[height=8cm,width=8cm]{parameter.eps}}}
%\end{picture}
%\end{center}
%\caption{\sl Parameters for a CR tetrahedron}\label{Figure:parameter}
%\end{figure}

We associate to each vertex $i$, in the edge $[ij]$, the invariant
$(ijkl)$ where the order of $k$ and $l$ is fixed by the right hand rule
 with the thumb
pointed from $j$ to $i$. A shortcut notation for the invariants is
therefore
$$
z_{ij}=(ijkl),
$$
the indices $kl$ being determined by the choice $ij$.

In order to give an explicit formula (see \cite{W2,Ge}),
denote by  $\hat p$ a lift of $p\in S^3$.
\begin{dfn} \label{polar}
Given two  points $p_i$ and $p_j$ in $S^3$,
a {\sl polar vector}
$c_{ij}$ is a  perpendicular vector to the complex line defined by
$p_i$ and $p_j$.
\end{dfn}
Observe that a polar vector is given by solving the equations
$\langle \hat p_i, c_{ij}\rangle=\langle \hat p_j, c_{ij}\rangle=0$ and
therefore
one can choose a polar vector to have coordinates in the field
defined by the conjugates of the coordinates of $\hat p_i$ and $\hat p_j$.

\begin{dfn}\label{polar}
 Consider a
generic configuration of four (ordered) points $(p_0,p_1,p_2,p_3)$
in $S^3$.  For $[p_i,p_j,p_k,p_l]\in \Delta O^+$,
 a positively oriented configuration, define
$$
z_{ij}=(ijlk)=\frac{\langle \hat p_l, c_{ij}\rangle \langle \hat p_k, \hat p_i \rangle}
{\langle \hat p_k, c_{ij}\rangle \langle \hat p_l, \hat p_i \rangle},
$$
\end{dfn}

It satisfies the relations $(ijlk)=(ijkl)^{-1}$ and
$$
z_{ij}z_{ji}=\overline {z_{kl}z_{lk}} \;.
$$

Moreover, we have the following description of the space of
configurations of four points.
\begin{prop} (cf. \cite{F2})
Configurations (up to translations by $PU(2,1)$) of four generic points
in $S^3$ are parametrised
by
$$
V\subset \left ( \C ^*\setminus \{ 1\}\right )^{12}
$$
with coordinates $z_{ij}$,  $1\leq i\neq j\leq 4$, defined by, for
$(i,j,k,l)$ an even permutation of $(0,1,2,3)$,
the usual similarity constraints
$$
z_{ik}=\frac{1}{1-z_{ij}}
$$
and the three complex equations
\begin{eqnarray}\label{eq:cr}
z_{ij}z_{ji}=\overline {z_{kl}z_{lk}}
\end{eqnarray}
with the exclusion of solutions such that
$z_{ij}z_{ji}z_{ik}z_{ki}z_{il}z_{li}=-1$ and $z_{ij}z_{ji}\in \R$.
\end{prop}

{\bf Remark}:
%\begin{enumerate}
%\item
The real solutions
parametrise
 configurations with
 four points
 contained in an $\R$-circle.
The solutions such that
$z_{ij}z_{ji}z_{ik}z_{ki}z_{il}z_{li}=-1$ and $z_{ij}z_{ji}\in \R$
correspond to
 degenerate hyperbolic ideal tetrahedra with  four points contained
in the boundary of a totally geodesic plane in real hyperbolic space.
They don't correspond to CR tetrahedra.
See \cite{Ge} for more details and \cite{F1, W1, W2, PP, PP1,FP} for other
descriptions.
%\end{enumerate}
\vspace{1cm}

We can also
describe generic configurations of four
points in $S^3$
by the following Lemma (see \cite[~Proposition 4.3]{F2}).
\begin{lem}\label{proposition:K}
One can always normalize a quadruple
$(p_0,p_1,p_2,p_3)$ of generic points in $S^3$ so that in the
compactified Heisenberg group they are
$$
p_0=\infty\ \ \ p_1=(0,0)\ \ \ p_2=(1,t) \ \ \ p_3=(z,s|z|^2),
$$
where $(z,t,s)$ are in the set
$$
 K=\{ (z,s,t)\in \C \times \R \times \R\ |\ z\ne 0,1\ \  \mbox{and}\ \ \overline{z}\frac{s+i}{t+i}\ne 1\ \}.
$$
In that case, the invariants of the configuration are
$$
 z_{01}=z,\; z_{10}=\frac{\overline{z}(s+i)}{t+i},\; z_{23}=\frac{z[(t+i)-\overline{z}(s+i)]}{(z-1)(t-i)},\;
 z_{32}=\frac{\overline{z}(z-1)(s-i)}{(t+i)-\overline{z}(s+i)}.
$$
\end{lem}
So the set $V$ in the above proposition is homeomorphic to the
set $K$ defined above.

\begin{dfn}\label{crinv}
To each configuration of four generic points $(p_0,p_1,p_2,p_3)$, we
define an element
$$
\beta(p_0,p_1,p_2,p_3):=
[z_{01}]+[z_{10}]+[z_{23}]+[z_{32}]
$$
in $\mathcal{P}(\C)$, the pre-Bloch group.
\end{dfn}

%\item The Korany-Reimann cross-ratio is
%$$
%KR(p_1,p_2,p_3,p_4)=(p_1,p_2,p_3,p_4)(p_2,p_1,p_4,p_3).
%$$

\section{A $\mathcal{P}(\C)$-valued CR invariant.}

In this section, let $G=PU(2,1)$.  For $n\geq 0$, we define
$C_n(S^3)$ to be the free abelian group
generated by the set of all generic $(n+1)$  ordered points in
$S^3$ (see Definition \ref{definition:generic}).

 The group $G$ acts on $S^3$ and therefore
 it acts diagonally on $C_n(S^3)$, which gives $C_n(S^3)$ a
left $G$-module structure.\\

We define the differential $d_n: C_n(S^3)\rightarrow C_{n-1}(S^3)$
by
$$
 d_n(p_0,\dots,p_n)=\sum_{i=0}^{n}(-1)^i (p_0,\dots,\hat{p_i},
 \dots, p_n),
$$
 then we can check that every $d_n$ is a $G$-module homomorphism and
$d_n\circ d_{n+1}=0$.
 Hence we have the $G$-complex
 $$
    C_{\bullet}(S^3): \cdots \rightarrow C_n(S^3)\rightarrow C_{n-1}(S^3)\rightarrow
    \cdots \rightarrow C_0(S^3).
 $$
We define the augmentation map $\epsilon: C_0(S^3)\rightarrow \Z $ by $\epsilon(x)=1$ for each $x\in S^3$.

\begin{lem} [cf. \cite{S1} pg. 221]
The augmentation complex $C_{\bullet}(S^3)\rightarrow \Z \rightarrow
0$ is exact.
\end{lem}

\Pf Let $C_{-1}(S^3)=\Z$. For $z=\sum_{i=1}^{m} n_i(p_0^{i},\dots,
p_k^{i})\in \text{ker} d_k$, we can choose a point $p\in S^3$ such
that if $k>0$, $p,p_0^{i},\dots, p_k^{i}$ are generic for all $i$;
if $k=0$, $p,p_0^{i}$ are distinct for all $i$. Hence
$(p,p_0^{i},\dots, p_k^{i})\in C_{k+1}(S^3)$ and we can check
directly that $d(\sum_{i=1}^{m} n_i(p,p_0^{i},\dots, p_k^{i}))=z$.
\EPf

For a left $G$-module $M$, we denote $M_G$ its group of co-invariants, that
is,
\[
 M_{G}=M/\langle gm-m, \forall g\in G, m\in M \rangle,
\]
where $\langle gm-m, \forall g\in G, m\in M \rangle$ is the
submodule of $M$ generated by all the elements of the form
$gm-m,g\in G, m\in M$. Take the co-invariants of the complex
$C_{\bullet}(S^3)$, we get the induced complex:
$$
    C_{\bullet}(S^3)_G: \cdots \rightarrow C_n(S^3)_G\rightarrow C_{n-1}(S^3)_G\rightarrow
    \cdots \rightarrow C_0(S^3)_G,
$$
with differential $\bar{d}_n: C_n(S^3)_G\rightarrow
C_{n-1}(S^3)_G$ induced by $d_n$. Since $G$ acts double
transitively on $S^3$, we see that $C_n(S^3)_G=\Z$ if $n\leq 1$, and
the differential $\bar{d}_1: C_1(S^3)_G\rightarrow
C_{0}(S^3)_G$ is zero. The equivalence class of
three generic points under the action of $G$ is determined by their
Cartan invariant (see Proposition
\ref{proposition:cartan}), so we get
$$
  C_2(S^3)_G=\Z[\R].
$$
Explicitly, a triple $(p_0,p_1,p_2)\in C_2(S^3)$ determines an element
$t=\tan \mathbb{A}(p_0, p_1, p_2)\in \Z[\R]$.  In normalized coordinates,
$(\infty,0,(1,t))\in C_2(S^3)$ is represented by $t\in \Z[\R]$.
 The differential
$$
\bar{d}_2: C_2(S^3)_G\rightarrow
C_{1}(S^3)_G=\Z
$$ is given on generators by $\bar{d}_2(p_0,p_1,p_2)=1$.\\

By Lemma \ref{proposition:K}, we can describe  $C_3(S^3)_G$
as follows:
$$
 C_3(S^3)_G= \Z[K],
$$
where
$$
 K=\{(z,s,t)\in \C \times \R \times \R|z\ne 0,1; \overline{z}\frac{s+i}{t+i}\ne 1\}.
$$
 In normalized coordinates,
$(\infty,0,(1,t),(z,s|z|^2))\in C_3(S^3)$ is represented by
$(z,s,t)\in \Z[K]$, so that
$$C_3(S^3)_G=
\bigoplus_{(z,s,t)\in K}\Z\cdot(\infty,0,(1,t),(z,s|z|^2)).
$$

By $\Z$-linear extension of Definition \ref{crinv} to ${C_3(S^3)_G}$,
we have a homomorphism
$$
  \bar \beta: {C_3(S^3)_G}\rightarrow \mathcal{P}(\C)
$$
given by
$$
 \sum_{l=1}^{k}n_l(p_0^l,p_1^l,p_2^l,p_3^l)\mapsto
\sum_{l=1}^{k}n_l\beta((p_0^l,p_1^l,p_2^l,p_3^l)).
$$

\begin{lem}
$\bar \beta(\text{Im} (\bar{d}_4))=0$ in $\mathcal{P}(\C)$.
\end{lem}

\Pf We need to show that the images are in the subgroup
generated by the 5-term relations. This is a special case of Theorem \ref{beta} and follows from Theorem
5.2 in \cite{F3}.
\EPf

 Therefore, $\bar \beta$ induces a well-defined homomorphism:
$$
 c: \frac{C_3(S^3)_G}{\text{Im}(\bar{d}_4)}\rightarrow
\mathcal{P}(\C)
$$
defined by
$$
c([ p_0,p_1,p_2,p_3])=\bar \beta(p_0,p_1,p_2,p_3),
$$
where $[ p_0,p_1,p_2,p_3]$ denotes the equivalence class of
$(p_0,p_1,p_2,p_3)$ in the quotient group. When we restrict
to $H_3(C_{\bullet}(S^3)_G)$, we get a homomorphism:
$$
 c: H_3(C_{\bullet}(S^3)_G)\rightarrow \mathcal{P}(\C).
$$

\begin{rmk}
In the real hyperbolic case, we have $G=PSL_2(\C)$ acting on $\C P^1$ and the corresponding $G$-complex. In that case,
we know that $H_3(C_{\bullet}(\C P^1)_{PSL_2(\C)})$ is equal to $\mathcal{P}(\C)$.
\end{rmk}

Normalizing the four points $p_0=\infty$, $p_1=0$, $p_2=(1,t)$,
$p_3=(z,s|z|^2)$ according to Lemma \ref{proposition:K}, we obtain
the following invariants;
$$
 z_{01}=z,\; z_{10}=\frac{\overline{z}(s+i)}{t+i},\; z_{23}=\frac{z[(t+i)-\overline{z}(s+i)]}{(z-1)(t-i)},\;
 z_{32}=\frac{\overline{z}(z-1)(s-i)}{(t+i)-\overline{z}(s+i)}
$$
and the homomorphism
$c:\frac{C_3(S^3)_G}{\text{Im}(\bar{d}_4)}\rightarrow
\mathcal{P}(\C)$ can be expressed as
$$
  c([ p_0,p_1,p_2,p_3])=[z_{01}]+[z_{10}]+[z_{23}]+[z_{32}].
$$

%By standard homological algebra, for every $p$, there is a canonical
%homomorphism:
%$$
%\alpha: H_p(G;\Z)\rightarrow H_p(C_{\bullet}(S^3)_G)=\frac{\text{Ker}(\overline{d_3})}{\text{Im}(\overline{d_4})}.
%$$

\begin{dfn} Given a configuration $[p_0,p_1,p_2,p_3]\in C_3(S^3)_G$,
we say it is
 symmetric if it lies in the kernel of
$\bar{d}_3: C_3(S^3)_G\rightarrow C_2(S^3)_G$.
\end{dfn}

Note that this definition is related to the definition of a
symmetric tetrahedron introduced in \cite[~Section 4.3]{F2}.
There, a configuration $[p_0,p_1,p_2,p_3]$ is symmetric if
there exists an anti-holomorphic involution $\phi$ such that
$\phi(p_i)=p_j$ and $\phi(p_k)=p_l$ for $\{i,j,k,l\}=\{0,1,2,3\}$.
Three possible symmetries may appear, namely
$\phi(p_0)=p_1$, $\phi(p_0)=p_3$ or $\phi(p_0)=p_3$.  They are characterized
by $\mathbb{A}(p_0,p_1,p_2)=\mathbb{A}(p_0,p_1,p_3)$,
 $\mathbb{A}(p_0,p_2,p_3)=\mathbb{A}(p_0,p_1,p_3)$ and
$\mathbb{A}(p_0,p_1,p_2)=-\mathbb{A}(p_0,p_2,p_3)$ respectively.
But the definition above concerns only the first two of them as shown in the next Lemma.

\begin{lem}\label{lem7.3}
$[p_0,p_1,p_2,p_3]\in C_3(S^3)_G$ is symmetric if and only if the
Cartan invariants satisfy
$$\mathbb{A}(p_0,p_1,p_2)=\mathbb{A}(p_0,p_1,p_3)\;\;
\text{or} \;\; \mathbb{A}(p_0,p_2,p_3)=\mathbb{A}(p_0,p_1,p_3).$$
\end{lem}

\Pf By the definition,
$$
 \bar{d}_3([p_0,p_1,p_2,p_3])=[p_1,p_2,p_3]-[p_0,p_2,p_3]+[p_0,p_1,p_3]-[p_0,p_1,p_2].
$$
Hence it lies in kernel of $\bar{d}_3$ if and only if in
$C_2(S^3)_G$, we have either
$$
  [p_1,p_2,p_3]=[p_0,p_2,p_3], \;\; \text{and} \;\; [p_0,p_1,p_3]=[p_0,p_1,p_2]
$$
or
$$
  [p_1,p_2,p_3]=[p_0,p_1,p_2], \;\; \text{and} \;\; [p_0,p_2,p_3]=[p_0,p_1,p_3].
$$
We know that in $C_2(S^3)_G$, two elements are equal if and only if
their Cartan invariants are the same. It is known that the Cartan
invariants satisfy the following cocycle conditions (\cite[~Page
219]{G}:
$$
\mathbb{A}(p_1,p_2,p_3)-\mathbb{A}(p_0,p_2,p_3)+\mathbb{A}(p_0,p_1,p_3)-\mathbb{A}(p_0,p_1,p_2)=0.
$$
Now the lemma follows. \EPf

\begin{prop}\label{prop7.4}
Given $[p_0,p_1,p_2,p_3]\in C_3(S^3)_G$ as above. Then it is
symmetric if and only if one of the following two equivalent
conditions holds:\\
(1). $t=s$ or $t+s-2(s\cdot\Re(z)+\Im(z))=0$;\\
(2). $|z_{01}|=|z_{32}|$.
\end{prop}

\Pf By \cite[~Proposition 4.6]{F2}, we have
$$
\tan\mathbb{A}(p_1,p_2,p_3)=\frac{2(s-t)\Re z+2(1+ts)\Im z+t(1+s^2)|z|^2-s(1+t^2)}{|(s-i)z+t-i|^2};
$$
and
$$
\tan\mathbb{A}(p_0,p_2,p_3)=\frac{|z|^2s-t+2\Im z}{|z-1|^2}, \; \tan\mathbb{A}(p_0,p_1,p_3)=s,\;
\tan\mathbb{A}(p_0,p_1,p_2)=t.
$$
By Lemma \ref{lem7.3}, we see that it is symmetric if and only if
$t=s$ or $s=\frac{|z|^2s-t+2\Im z}{|z-1|^2}$, which is condition
(1). By the definition of $z_{01}$ and $z_{32}$, a direct
computation shows that $|z_{01}|^2=|z_{32}|^2$ if and only if
$(t-s)(t+s-2(s\cdot\Re(z)+\Im(z)))=0$. \EPf

Recall that
$$
 \mathcal{P}(\C)=\mathcal{P}(\C)^{+}+\mathcal{P}(\C)^{-}
$$
is the decomposition of $\mathcal{P}(\C)$ into the two subgroups preserved by
the complex conjugation involution.

 \begin{prop}
(1). The image of the homomorphism $c:H_3(C_{\bullet}(S^3)_G)\rightarrow
\mathcal{P}(\C)$ contains $\mathcal{P}(\C)^{+}$, the invariant subgroup
of $\mathcal{P}(\C)$ under complex conjugation.\\
(2). Suppose $[p_0,p_1,p_2,p_3]\in C_3(S^3)_G$ is symmetric. Then its
image under the homomorphism $c$ lies in  $\mathcal{P}(\C)^{+}$.
\end{prop}

\Pf For (1), we show that the subgroup generated by the images of the
symmetric elements contains $\mathcal{P}(\C)^{+}$. Consider the
elements $[p_0,p_1,p_2,p_3]$ of the form
$$
 p_0=\infty, p_1=0, p_2=(1,t), p_3=(z,t|z|^2); \; z\in \C-\{0,1\}, t\in \R.
$$
By Proposition \ref{prop7.4}, they are symmetric. By
\cite[~Corollary 4.11]{F2}, we find the invariants:
$$
z_{01}=z, z_{10}=\overline{z}, z_{23}=ze^{i \theta}, z_{32}=\overline{z}e^{-i \theta};\;\; z\in \C-\{0,1\}, \theta\in
\R.
$$
Therefore, we have
\begin{equation}\label{e1}
 c([p_0,p_1,p_2,p_3])=[z]+[\overline{z}]+[ze^{i\theta}]+[\overline{z}e^{-i\theta}]
\end{equation}

By \cite[~Theorem 4.16]{Sah}, $\mathcal{P}(\C)$ is a $\Q$-vector
space.
Let $B=\sum_{i=1}^{k}n_i[a_i]\in
\mathcal{P}(\C)^{+}$. Then
$\sigma(B)=\sum_{i=1}^{k}n_i[\overline{a_i}]=B$. Hence
$$
 B=\frac{1}{2}(B+\sigma(B))=\frac{1}{2}\sum_{i=1}^{k}n_i([a_i]+[\overline{a_i}]).
$$
Choose $b_i\in \C$ such that $a_i=b_{i}^2$ and therefore
$\overline{a_i}=\overline{b_i}^2$. By \cite[~Theorem 5.23]{DS}, we
know that in $\mathcal{P}(\C)$, $[a^2]=2([a]+[-a])$. Therefore,
$$
 B=\sum_{i=1}^{k}n_i([b_i]+[\overline{b_i}]+[-b_i]+[-\overline{b_i}]).
$$
Now in (\ref{e1}), if we choose $z=b_i$ and $e^{i\theta}=-1$, we see
that the first part follows.

For (2), let $T=c([p_0,p_1,p_2,p_3])=[z_{01}]+[z_{10}]+[z_{23}]+[z_{32}]$. Using the 5-term relations, a direct
computation shows that in $\mathcal {P}(\C)$, the difference $\sigma(T)-T$ is equal to
$$
 [\frac{z_{32}(1-z_{10})}{\overline{z}_{01}(1-\overline{z}_{23})}]+
          [\frac{\overline{z}_{32}(1-\overline{z}_{23})}{z_{01}(1-z_{10})}]+
          [\frac{z_{01}(1-z_{23})}{\overline{z}_{32}(1-\overline{z}_{10})}]+
          [\frac{\overline{z}_{01}(1-\overline{z}_{10})}{z_{32}(1-z_{23})}]
$$
Put
$$
 a=\frac{z_{32}(1-z_{10})}{\overline{z}_{01}(1-\overline{z}_{23})},
 b=\frac{\overline{z}_{32}(1-\overline{z}_{23})}{z_{01}(1-z_{10})},
$$
we can rewrite
$$
  \sigma(T)-T=[a]+[b]+[(\overline{a})^{-1}]+[(\overline{b})^{-1}],\;\text{and} \;ab=\frac{|z_{32}|^2}{|z_{01}|^2}.
$$
Since $[p_0,p_1,p_2,p_3]$ is symmetric, by Proposition \ref{prop7.4}, $ab=1$, i.e. $b=a^{-1}$. On the other hand, we
know in $\mathcal {P}(\C)$,
$$
  [z]+[z^{-1}]=0,
$$
therefore
$$
  \sigma(T)-T=[a]+[a^{-1}]+[(\overline{a})^{-1}]+[\overline{a}]=0.
$$
That is, $\sigma(T)=T$ and the image lies in $\mathcal{P}(\C)^{+}$.
\EPf

\begin{rmk}
The previous proposition shows that the homomorphism
$c:H_3(C_{\bullet}(S^3)_G)\rightarrow \mathcal{P}(\C)$ is
non-trivial, and $\mathcal{P}(\C)^{+}$ is equal to the image of the subgroup generated by the symmetric
configurations. It would be very interesting to determine its kernel and image.
\end{rmk}

Let $M$ be a $3$-dimensional spherical CR manifold (possibly non-compact).
Suppose it has a
triangulation
consisting of finitely many CR tetrahedra, say
$M=\Delta_1\cup \cdots \cup \Delta_k$ with each $\Delta_i$ a CR
tetrahedron (here we suppose that the triangulation is ideal if the
manifold is non-compact). Denote $p_0^i,p_1^i,p_2^i,p_3^i$ the four vertices of
$\Delta_i$ where the order of the vertices is consistent with the
orientation.
\begin{dfn}
Let $M$ be as above with the CR triangulation. Define an element
$[M]\in \frac{C_3(S^3)_G}{\text{Im}(\bar{d}_4)}$ by
$$
[M]:=\sum_{i=1}^{k}[p_0^i,p_1^i,p_2^i,p_3^i].
$$
\end{dfn}
\begin{lem}
1. $[M]$ is independent of the triangulation.\\
2. $[M]\in H_3(C_{\bullet}(S^3)_G)$, that is,
$\bar{d}_3([M])=0$.
\end{lem}

\Pf 1. If $M$ is closed, we know that two different triangulations
can be obtained from one to the other by Pachner moves. If $M$ is
not closed, by Theorem \ref{tr}, the same result holds. From the
definition, it is clear one Pachner move gives an element of
$\text{Im}(\bar{d}_4)$. Hence $[M]$ is independent of the
triangulation.

2. Since $M$ is triangulated, their faces are
matched and the terms in ${\bar{d}_3}([M])$ are canceled out in
pairs. \EPf

In the following definition we suppose that a triangulation of a non-compact
manifold is ideal.

\begin{dfn}\label{def311}
 Let $M$ be a $3$-dimensional spherical CR manifold with a CR
triangulation. We define $\beta(M):=c([M])\in \mathcal{P}(\C)$.
\end{dfn}

Now recall the Bloch-Wigner function (See Definition \ref{BWf})
$ D: \mathcal{P}(\C)\rightarrow \R$.

\begin{thm}\label{vol}
Let $M$ be as in the definition \ref{def311}. Then $D(\beta(M))=0$.
\end{thm}

\Pf
For a configuration $[p_0,p_1,p_2,p_3]\in C_3(S^3)_G$ with its cross-ratios $z_{ij}$, by remark $3$ of section $6$ in
\cite[~page 14]{F3}, we have the following identity:
$$
-e^{2i\mathbb{A}(p_i, p_j, p_k)}=z_{il}z_{jl}z_{kl},
$$
where $[i,j,k,l]$ is an even permutation of $[0,1,2,3]$.
By \cite[~Proposition 6.5]{F3}, we obtain:
$$
 2D(c([p_0,p_1,p_2,p_3]))=D(-e^{2i\mathbb{A}(p_1,p_2,p_3)})+
D(-e^{2i\mathbb{A}(p_0,p_3,p_2)})+
 D(-e^{2i\mathbb{A}(p_0,p_1,p_3)})+D(-e^{2i\mathbb{A}(p_0,p_2,p_1)}).
$$
Let $M=\Delta_1\cup \cdots \cup \Delta_k$ with each $\Delta_i$ a CR
tetrahedron. Denote $p_0^i,p_1^i,p_2^i,p_3^i$ the four vertices of
$\Delta_i$ and the order of the vertices is consistent with the
orientation. Then by definition
$$
D(\beta(M))=\sum_{i=1}^{k}D(c([p_0^i,p_1^i,p_2^i,p_3^i])).
$$
Since $M$ is glued by the tetrahedra $\Delta_i$, their faces are glued in pairs, and the corresponding Cartan
invariants are equal. By the above formula of $2D(c([p_0,p_1,p_2,p_3]))$, we see that the terms in $D(\beta(M))$ are
canceled in pairs. Hence it is zero.
\EPf
\begin{rmk}
Note that $D(\beta(M))$ is the CR volume of $M$ defined in \cite{F3}. The above theorem says that it is always zero.
This is exactly the opposite to the real hyperbolic case where the volume is never zero.
\end{rmk}

\section{An invariant in  $\mathcal{B}(k)$}

In this section, we will discuss when the $\mathcal{P}(\C)$-valued invariant
 defined above can be defined
 in the Bloch group $\mathcal{B}(k)$ for $k$ a number field. We first discuss the real
hyperbolic case, where it is known that the invariant always lies in $\mathcal{B}(\C)$.

Observe first that given a cross-ratio structure $(T, \bf{X})$
 associated to a triangulation we may associate a field
$k_{_{\bf{X}}}=\Q(z_1, \cdots )$, where $z_1,\cdots$
are all the cross-ratios. It is clearly preserved by taking finite coverings
of the structure.
From Proposition 4.2 in \cite{F3} we obtain:

\begin{prop}
The field  $k_{_{\bf{X}}}$ is invariant under elementary moves.
\end{prop}

For the case of an ideal real hyperbolic triangulation see
\cite{NR}. In particular, one can compare $k_{_{\bf{X}}}$ to a holonomy
representation defined by an ideal triangulation.  Recall that the invariant
trace field of a representation $\rho : \Gamma \rightarrow PSL(2, \C)$ is
given by taking lifts $\tilde g\in SL(2,\C)$ of elements $g\in PSL(2, \C)$ and
defining (\cite{R})
$$
k_\rho=\Q\left(\left\{\  Tr(\tilde g^2)\ |\ g\in \Gamma\ \right\}\right).
$$
For an ideal triangulation of a finite volume cusped hyperbolic manifold, the
field obtained by adjoining the
cross-ratios of the ideal tetrahedra and the field obtained from a
holonomy representation $\rho$ are the same, that
is
$k_\rho = k_{_{\bf{X}}}$ (\cite{NR} Theorem 2.4 pg. 277). Moreover, if we choose
one tetrahedron with one of its faces normalized to be,
in homogeneous coordinates of $\C P^1$, $[1,0], [0,1], [1,1]$
(that is, $\infty, 0, 1$ in $\C\cup \{\infty\}$), then the coordinates
$[z_i,1]$ of the
vertices of the other tetrahedra obtained by developing the triangulation
are all in the field  $k_{_{\bf{X}}}$ and any side pairing $g\in PSL(2,\C)$
which identifies two sides of the triangulation has a lift
with coefficients in the same field (see \cite{NR} lemma 2.5 pg. 278).

In the case of CR structures analogous results were proven in \cite{Ge1}.
Denote by $\tilde g\in SU(2,1)$ a lift of an element $g\in PU(2, 1)$.
Let the invariant
trace field of a representation $\rho : \Gamma \rightarrow PU(2,1)$ be defined
as (see \cite{Ge1, Mc})
$$
k_\rho=\Q\left(\left\{\  Tr((\tilde \rho(\gamma))^3)\ |\ \gamma\in \Gamma\ \right\}\right).
$$
For an ideal triangulation of a CR structure, we can
develop it in $S^3$ by choosing
one tetrahedron with one of its faces normalized to be,
in homogeneous coordinates of $\C P^2$,
$$[1,0,0], [0,0,1], [(-1+it)/2,1+it,1]
$$
(that is, $\infty, 0, 1+it$ in the compactified Heisenberg space
$\overline{\Heis}$), then the coordinates
$[z^1_i,z^2_i,1]$ of the
vertices of the other tetrahedra obtained by developing the triangulation
are all in the field  $k_{_{\bf{X}}}$.  We will call such a construction
a {\sl normalized development}.

\begin{prop}[\cite{Ge}]\label{coordinates}
If  $(T,{\bf X})$ is a CR triangulation then
$$
k_{_{\bf{X}}}=\Q(z^1_{i},z^2_{i},{\bar z}^1_{i},{\bar z}^2_{i}),
$$
where $[z^1_i,z^2_i,1]$ are as above, the coordinates of the vertices of a normalized development.  In particular,
$k_{_{\bf{X}}}$ is invariant under complex conjugation.
\end{prop}

%\begin{prop} [\cite{Ge} proposition 3.3]
%For an ideal triangulation of a CR structure with a Zariski dense holonomy representation such that each vertex of the developed triangulation is a parabolic fixed point, the field obtained by adjoining the cross-ratios of the ideal tetrahedra and the field obtained from the holonomy representation $\rho$ are the same, that is $k_\rho = k_{_{\bf{X}}}$.\end{prop}

%Moreover any side pairing $g\in PU(2,1)$ which identifies two sides of the triangulation has a lift $\tilde g$ with $(\tilde g)^3\in SU(2,1,k_{_{\bf{X}}})$ (see Lemma 3.4 in \cite{Ge1}).

%\begin{thm}[\cite{Ge1}]\label{Zariski} Let $\Gamma$ be a subgroup of $SU(2,1)$.  Suppose that $\Gamma$ contains a parabolic element.  If $\Gamma$ is Zariski dense then it can be conjugated to a subgroup of $SU(2,1)$ with coefficientsin the field $k_{_\Gamma}$.\end{thm} Here, $k_{_\Gamma}$ is the field defined as above by the identity representation.

\begin{dfn}\label{invariant} Let $(T, \bf{X})$ be a cross-ratio structure. Define
$$\beta_k (T, {\bf{X}})=
\sum_{s} ([z^s_{01}]+[z^s_{10}]+[z^s_{23}]+[z^s_{32}]
%+[z^s_{02}]+[z^s_{20}]+[z^s_{13}]+[z^s_{31}]+[z^s_{03}]+[z^s_{30}]+[z^s_{12}]+[z^s_{21}]
)\in
\mathcal{P}(k_{_{\bf{X}}}).
$$
%$$=\sum_s \sum_{i\neq j}[z^s_{ij}]$$
\end{dfn}

{\bf Remark}:\label{Remark}
This definition depends on a choice of edge in each simplex of a triangulation.
There are several ways to make it independent of the choice.  We can
use Lemma \ref{2beta} and define it to be $2\beta (T, {\bf{X}})$ in
$\mathcal{P} (k_{_{\bf{X}}}(\omega))$, where $\omega$ is a cubic root of unity.
 Or define it to be
$$\sum_{s} ([z^s_{01}]+[z^s_{10}]+[z^s_{23}]+[z^s_{32}]
+[z^s_{02}]+[z^s_{20}]+[z^s_{13}]+[z^s_{31}]+[z^s_{03}]+[z^s_{30}]+[z^s_{12}]+[z^s_{21}]
)\in
\mathcal{P}(k_{_{\bf{X}}}).
$$
Another definition that is independent of the particular choice
of edges in each tetrahedron starts  by defining
 a slightly smaller pre-Bloch group $\mathcal{P'}(k_{_{\bf{X}}})$
as in (\cite{NY} Definition 2.3 pg. 4) to be
 the quotient of
the free abelian group $\Z[k_{_{\bf{X}}}\cup \infty]$
by the subgroup generated by
5-term relations
and
$$
[1]=[0]=[\infty]=0.
$$
Both definitions of the pre-Bloch group will differ by a torsion subgroup.
With the definition above we have $[z]=[\frac{1}{1-z}]$ and therefore one can define $\sum_{s}
([z^s_{01}]+[z^s_{10}]+[z^s_{23}]+[z^s_{32}])$
as an invariant which does not depend on the choice of a pair of edges in each tetrahedron.

\subsection{The Real Hyperbolic invariant}

We will show for a finite
volume non-compact real hyperbolic 3-manifold with ideal triangulation,
the corresponding invariant lies in $B(\C)$. This follows from \cite{NZ}
(see a proof in \cite{NY}).  But we give
an elementary geometric proof of the identity which will give an idea of the
proof in the CR case.

Consider an ideal triangulation of a 3-manifold.
We order the simplices by choosing arbitrarily a first one and then
a second one with a common face.
Having chosen $n$ simplices we choose the $(n+1)$-th as a simplex with one face in common
with the union of the previous ones or, if there are two common faces, they share
a common edge (if there are three common faces, the edges defined by each pair should
be common too).  In this way we obtain a polyhedron homeomorphic to a 3-ball
with face pairings on its boundary (homeomorphic to a sphere).

Suppose the triangulation has a hyperbolic cross-ratio structure.
To each order on the simplices as above,
we can associate a well defined map from the 0-skeleton of each simplex to $\C P^1$.
We can associate with the first simplex the points $\infty, 0, 1, z$ and
proceed determining the other vertices which are clearly defined by the cross-ratios.
In fact, the map is determined by the chosen order on the simplices
 and the initial map defined on the 0-skeleton of the
 first simplex.

\begin{dfn}
The holonomy group of the hyperbolic structure is
the group generated by the face
pairing transformations
of the polyhedron.
\end{dfn}

%One can always lift the representation holonomy
% of a hyperbolic 3-manifold
%into $PSL(2,\C)$  to a representation into $SL(2,\C)$ (see \cite{CS,Cu}).

Let $p_i\in \C P^1$ be the vertices of the development map as above.
Choose a lift $v_i\in \C^2$ for each vertex $p_i$. The vertices
of the faces of the polyhedron
are identified by side pairings in $PSL(2,\C)$ but their lifts to $\C^2$
might not be
 identified by lifts of the side pairings to $SL(2,\C)$.

\begin{dfn}
We call a lift of the vertices of a development of a
 triangulation {\sl special}
if there are  lifts of the side pairings which preserve the lifted vertices
 up to multiplication by $-1$.
\end{dfn}

\begin{lem}
Any ideal triangulation of a finite volume 3-manifold has a special lift.
\end{lem}
\Pf We consider a finite ideal triangulation of the 3-manifold and its
development by ideal tetrahedra in hyperbolic space.  We are only  interested
in its vertices in $\C P^1$.  Let $\Gamma\in PSL(2,\C)$ be the holonomy group
of the hyperbolic structure.

%We write $\hat g\in SL(2,\C)$ for the lift
%of an  element $g\in PSL(2,\C)$.

Choose one vertex $p\in \C P^1$  and  a lift $v\in \C^2$ of $p$.
 Without loss of generality,
 we suppose that $v=(1,0)$ ($p=\infty$).  The elements of $SL(2,\C)$
 which are lifts of elements of the parabolic group fixing $p$ are of the form
$$
 \pm \left [ \begin{array}{cc}
                      1 & \star\\
                      0 & 1
                   \end{array}
           \right ].
$$
Consider all vertices
 identified to $p$ by the holonomy group.
That is, the
other vertices $p_i$
are obtained by $p_i=g_ip$ for $g_i\in \Gamma$.
  Let
$\hat g_i\in SL(2,\C)$ be a lift of $g_i$ and let $v_i=\hat g_iv$.
We first observe that the definition is compatible up to a multiplication by
$-1$.  Indeed, if $g_kp=g_lp$ then $g_l^{-1}g_kp=p$ and therefore
$g_l^{-1}g_k$ is parabolic (or the identity) and $\hat g_kv=\pm\hat g_lv$.

Let $g\in  \Gamma$ be a side pairing of the polyhedra obtained by
the development map.  We have
$
gg_ip=g_jp,
$
therefore $(g_j)^{-1}gg_ip=p$ and we conclude that $(g_j)^{-1}gg_i$ is parabolic
and then its lift $(\hat g_j)^{-1}\hat g\hat g_i$ is a parabolic element (or the identity element) fixing $v$ up to
sign. This implies that $\hat g\hat g_iv =\pm \hat g_j v$ and we conclude that
$$
\hat g v_i=\pm v_j.
$$
We do that for each cycle of vertices and
obtain that the lift is special.
\EPf

Here we use the following observation in \cite{DZ}.
 Let $v_i=(v_i^1,v_i^2)\in \C^2$ for $0\leq i\leq 3$ and define the determinant
$$
[v_i,v_j]= \left | \begin{array}{cc}
                      v_i^1 & v_j^1\\
                      v_i^2 & v_j^2
                   \end{array}
           \right |.
$$
Suppose $(v_0,v_1,v_2,v_3)$ is a quadruple of points in $\C^2$ so that
 they define a quadruple of pairwise distinct points
 $(p_0,p_1,p_2,p_3)$ in $\C P^1$, where $h(v_i)=p_i$ is the projection
in projective space.
Observe that
$$
[p_0,p_1,p_2,p_3]=\frac{[v_2,v_0][v_3,v_1]}{[v_3,v_0][v_2,v_1]}
$$
is the cross ratio of the projection of the four distinct points
in $\C P^1$.

%Let $p_1,p_2\in \C P^1$ be two distinct points. Let $(q_1,q_2)= g (p_1,p_2)$ for $g\in PSL(2,\C)$.
%Then, for any lift $v_i\in \C^2$ of $p_i$ and  $w_i\in \C^2$ of $q_i$ we
%have $[v_1,v_2]=\lambda [w_1,w_2]$ where $\lambda\in \C^*$ does not depend on $g$, only on the lifts.

If we choose a special hyperbolic lift of a triangulation,
 we obtain a well defined function (up to a sign) on the 1-simplices
of the triangulation.  That function recuperates
the cross ratio of the 3-simplices according
to the formula above.

\subsection{Bloch identity}

\begin{dfn}
Given a finite hyperbolic triangulation $T=(T_i, z_i)$ by ideal tetrahedra,
 we define its Bloch sum by
$$
\delta(\beta(T))= \sum_i z_i\wedge (1-z_i).
$$
\end{dfn}
We prove

\begin{thm}For a finite ideal triangulation $T$
of a finite volume hyperbolic manifold
$$
\delta(\beta(T))=0.
$$
\end{thm}
\Pf
Consider a development of the triangulation with vertices in $\C P^1$.  There
exists a special lift by the Proposition above.
We compute the following sum
$$
\sum_i [p^i_0,p^i_1,p^i_2,p^i_3]\wedge [p^i_0,p^i_2,p^i_3,p^i_1].
$$
For each tetrahedron $(p_0,p_1,p_2,p_3)$, let $v_i$ be the lift of $p_i$.  Then
$$
[p_0,p_1,p_2,p_3]=\frac{[v_2,v_0][v_3,v_1]}{[v_3,v_0][v_2,v_1]}.
$$
Therefore
$$
[p_0,p_1,p_2,p_3]\wedge [p_0,p_2,p_3,p_1]=\frac{[v_2,v_0][v_3,v_1]}{[v_3,v_0][v_2,v_1]}\wedge
\frac{[v_3,v_0][v_1,v_2]}{[v_1,v_0][v_3,v_2]}
$$
$$
= [v_0,v_2]\wedge [v_0,v_3] +[v_0,v_2]\wedge [v_2,v_1] -[v_0,v_2]\wedge [v_2,v_3]-[v_0,v_2]\wedge [v_0,v_1]
$$
$$
+\, [v_1,v_3]\wedge [v_0,v_3] +[v_1,v_3]\wedge [v_2,v_1] -[v_1,v_3]\wedge [v_2,v_3]-[v_1,v_3]\wedge [v_0,v_1]
$$
$$
-\, [v_1,v_2]\wedge [v_0,v_3] -[v_1,v_2]\wedge [v_2,v_1] +[v_1,v_2]\wedge [v_2,v_3]+[v_1,v_2]\wedge [v_0,v_1]
$$
$$
-\, [v_0,v_3]\wedge [v_0,v_3] -[v_0,v_3]\wedge [v_2,v_1] +[v_0,v_3]\wedge [v_2,v_3]+[v_0,v_3]\wedge [v_0,v_1]
$$
Observe that $(-a)\wedge b = a\wedge b$ and $[v,w]=-[w,v]$.  Therefore  we can write
$$
[p_0,p_1,p_2,p_3]\wedge [p_0,p_2,p_3,p_1] =
$$
$$
\ [v_0,v_1]\wedge [v_0,v_2] + [v_0,v_2]\wedge [v_0,v_3] +[v_0,v_3]\wedge [v_0,v_1]+
$$
$$
\  [v_1,v_2]\wedge [v_1,v_0] + [v_1,v_3]\wedge [v_1,v_2] +[v_1,v_0]\wedge [v_1,v_3]+
$$
$$
\  [v_2,v_0]\wedge [v_2,v_1] +[v_2,v_1]\wedge [v_2,v_3] +[v_2,v_3]\wedge [v_2,v_0]+
$$
$$
 [v_3,v_0]\wedge [v_3,v_2] +[v_3,v_2]\wedge [v_3,v_1] +[v_3,v_1]\wedge [v_3,v_0]
$$
This shows that each face of a fixed tetrahedron appears three times in the sum
but the terms corresponding to a common face between two tetrahedra appear
with opposite sign.  When the faces are identified by a side pairing,
the lift of the side pairing, at most, changes the sign of the lifted vertices.
But the terms of the sum above are invariant under sign change.
\EPf

\subsection{CR invariant}
Let $v_1,v_2\in \C^{2,1}$ be two vectors in the Hermitian space $\C^{2,1}$ with Hermitian
product $\langle \cdot , \cdot \rangle$.  We will write $v_{12}=v_1\boxtimes v_2$ for the
(alternating bilinear) Hermitian cross-product as defined in \ref{boxproduct}.
  It satisfies, for any $v_1,v_2,v_3\in \C^{2,1}$,
$$
\langle v_1, v_{23}\rangle = \langle v_2, v_{31}\rangle .
$$

For four generic points $(p_0,p_1,p_2,p_3)$ in $S^3$
and $v_i\in \C^{2,1}$ chosen lifts, put
$$
[p_0,p_1,p_2,p_3]= \frac{\langle v_3,v_{01}\rangle \langle v_2,v_{0}\rangle }
{\langle v_2,v_{01}\rangle \langle v_3,v_{0}\rangle }.
$$
Then by Definition \ref{crinv}, \ref{polar}, we have
$$
\beta(p_0,p_1,p_2,p_3)=
[p_0,p_1,p_2,p_3]+[p_1,p_0,p_3,p_2]+[p_2,p_3,p_0,p_1]+[p_3,p_2,p_1,p_0].
$$

The next proposition computes the image of $\beta(p_0,p_1,p_2,p_3)$ under the map
$$
\delta : \mathcal{P}(\C) \rightarrow \C^*\wedge \C^*
$$
defined on the generators of $\mathcal{P}(\C)$  by
$\delta([z])=z\wedge(1-z)$.  The main objective is to obtain an expression
which depends on  four oriented surface terms.  These terms will cancel
out when two configurations have a common face.

\begin{lem} Let $p_i\in S^3$, $0\leq i\leq 3$, be four
 pairwise distinct points
and $v_i\in \C^{2,1}$ chosen lifts (to simplify notations we will denote
$v_i$ simply $i$).  We have
$$
-\delta(\beta(p_0,p_1,p_2,p_3))=
\langle 3, 01\rangle\wedge \frac{\langle 0,1\rangle\langle 1, 3\rangle \langle 3, 0\rangle}{\langle 0,3\rangle\langle
3, 1\rangle \langle 1, 0\rangle}
+\langle 2, 01\rangle\wedge \frac{\langle 0,2\rangle\langle 2, 1\rangle \langle 1, 0\rangle}{\langle 0,1\rangle\langle
1, 2\rangle \langle 2, 0\rangle}
$$
$$
+\langle 3, 02\rangle\wedge \frac{\langle 0,3\rangle\langle 3, 2\rangle \langle 2, 0\rangle}{\langle 0,2\rangle\langle
2, 3\rangle \langle 3, 0\rangle}
+\langle 2, 31\rangle\wedge \frac{\langle 1,2\rangle\langle 2, 3\rangle \langle 3, 1\rangle}{\langle 1,3\rangle\langle
3, 2\rangle \langle 2, 1\rangle}
$$
$$
+\langle 2,0\rangle\wedge \langle 3,0\rangle+
\langle 1, 0\rangle\wedge \langle 2, 0\rangle
+\langle 3,0\rangle\wedge \langle 1,0\rangle
$$
$$
+\langle 3,1\rangle\wedge \langle 2,1\rangle+
\langle 0, 1\rangle\wedge \langle 3, 1\rangle
+\langle 2,1\rangle\wedge \langle 0,1 \rangle
$$
$$
+\langle 0,2\rangle\wedge \langle 1,2\rangle+
\langle 3, 2\rangle\wedge \langle 0, 2\rangle
+\langle 1,2\rangle\wedge \langle 3,2 \rangle
$$
$$
+\langle 1,3\rangle\wedge \langle 0,3\rangle+
\langle 2, 3\rangle\wedge \langle 1, 3\rangle
+\langle 0,3\rangle\wedge \langle 2,3 \rangle.
$$

\end{lem}
\Pf
To simplify notation we write
$$
[p_0,p_1,p_2,p_3]=[0,1,2,3]=\frac{\langle 3,{01}\rangle \langle 2,{0}\rangle }
{\langle 2,{01}\rangle \langle 3,{0}\rangle }.
$$
Also, observe that
$$
-\delta([0,1,2,3])=[0,1,2,3]\wedge [0,2,3,1].
$$
Therefore we need to compute
$$
-\delta(\beta(p_0,p_1,p_2,p_3))=[0,1,2, 3]\wedge [ 0, 2, 3, 1]+
[ 1, 0, 3, 2]\wedge [ 1, 3, 2, 0]
$$
$$
+[ 2, 3, 0, 1]\wedge [ 2, 0, 1, 3]+
[ 3, 2, 1, 0]\wedge [ 3, 1, 0, 2]
$$
$$
=
\frac{\langle 3, 01\rangle\langle 2, 0\rangle}{\langle 2, 01\rangle\langle 3, 0\rangle}
\wedge
\frac{\langle 1, 02\rangle\langle 3, 0\rangle}{\langle 3, 02\rangle\langle 1, 0\rangle}
+\frac{\langle 2, 10\rangle\langle 3, 1\rangle}{\langle 3, 10\rangle\langle 2, 1\rangle}
\wedge
\frac{\langle 0, 13\rangle\langle 2, 1\rangle}{\langle 2, 13\rangle\langle 0, 1\rangle}
$$
$$
+\frac{\langle 1, 23\rangle\langle 0, 2\rangle}{\langle 0, 23\rangle\langle 1, 2\rangle}
\wedge
\frac{\langle 3, 20\rangle\langle 1, 2\rangle}{\langle 1, 20\rangle\langle 3, 2\rangle}
+\frac{\langle 0, 32\rangle\langle 1, 3\rangle}{\langle 1, 32\rangle\langle 0, 3\rangle}
\wedge
\frac{\langle 2, 31\rangle\langle 0, 3\rangle}{\langle 0, 31\rangle\langle 2, 3\rangle}
$$
We may use the distributive property of the wedge product to obtain a sum
of three types of terms.  The first type one is of the form
$\langle i,jk \rangle\wedge \langle i',j'k'\rangle$.  The second one is of the form
$\langle i,jk \rangle\wedge \langle n,m\rangle$.  The third type is
$\langle i,j \rangle\wedge \langle n,m\rangle$.

Using the property $\langle 1, {23}\rangle = \langle 2, {31}\rangle=
-\langle 1, {32}\rangle$ and the fact that $(-1)\wedge z=0$ for all $z$ we can treat
$\langle i,jk \rangle$ as invariant under all permutations.

We obtain by a straightforward computation that all terms of type 1 cancel out.
Also, the terms of type 2 can be written in a more concise form as
$$
\langle 3, 01\rangle\wedge \frac{\langle 0,1\rangle\langle 1, 3\rangle \langle 3, 0\rangle}{\langle 0,3\rangle\langle
3, 1\rangle \langle 1, 0\rangle}
+\langle 2, 01\rangle\wedge \frac{\langle 0,2\rangle\langle 2, 1\rangle \langle 1, 0\rangle}{\langle 0,1\rangle\langle
1, 2\rangle \langle 2, 0\rangle}
+\langle 3, 02\rangle\wedge \frac{\langle 0,3\rangle\langle 3, 2\rangle \langle 2, 0\rangle}{\langle 0,2\rangle\langle
2, 3\rangle \langle 3, 0\rangle}
+\langle 2, 31\rangle\wedge \frac{\langle 1,2\rangle\langle 2, 3\rangle \langle 3, 1\rangle}{\langle 1,3\rangle\langle
3, 2\rangle \langle 2, 1\rangle}
$$

Finally, the terms of type 3 are the following
$$
\langle 2,0\rangle\wedge \langle 3,0\rangle+
\langle 1, 0\rangle\wedge \langle 2, 0\rangle
+\langle 3,0\rangle\wedge \langle 1,0\rangle
$$
$$
+\langle 3,1\rangle\wedge \langle 2,1\rangle+
\langle 0, 1\rangle\wedge \langle 3, 1\rangle
+\langle 2,1\rangle\wedge \langle 0,1 \rangle
$$
$$
+\langle 0,2\rangle\wedge \langle 1,2\rangle+
\langle 3, 2\rangle\wedge \langle 0, 2\rangle
+\langle 1,2\rangle\wedge \langle 3,2 \rangle
$$
$$
+\langle 1,3\rangle\wedge \langle 0,3\rangle+
\langle 2, 3\rangle\wedge \langle 1, 3\rangle
+\langle 0,3\rangle\wedge \langle 2,3 \rangle.
$$
\EPf

Consider a CR triangulation and Let $p_i\in
S^3$ be the vertices of the developement map as in the real
hyperbolic case. Choose a lift $v_i\in \C^3$ for each vertex $p_i$.
The vertices of the faces of the polyhedron are identified by side
pairings in $PU(2,1)$.

\begin{dfn}
Define the  holonomy group of a CR triangulation to be
the group generated by the face
pairing transformations
of the polyhedron.
\end{dfn}

\begin{dfn}
We call a lift of the vertices of a development of a
CR triangulation {\sl special }
if there are  lifts of the side pairings preserving the lifted vertices
 up to multiplication by a root of unity.
\end{dfn}

\begin{lem}\label{special}
If a CR triangulation is such that
each vertex of the developed triangulation
is the fixed point of purely parabolic
elements (or the identity) of the holonomy group
then it has a special lift.
\end{lem}

\Pf The proof is the same as in the case of real hyperbolic
structures.  That is, we divide the vertices into classes, each one
obtained by translating a fixed vertex by elements of the holonomy
group.
That is, fixing a vertex $p\in S^3$ we consider all vertices
 identified to $p$ by the holonomy group.
That is, the
other vertices $p_i$
are obtained by $p_i=g_ip$ for $g_i\in \Gamma$. In order to
verify the compatibility, suppose
$g_ip=g_jp$. Then $g_j^{-1}g_ip=p$ and by hypothesis $g_j^{-1}g_i$
is purely parabolic. Without loss of generality, we suppose that $p=\infty$ and
let
$\hat g_i\in SU(2,1)$ be a lift of $g_i$ and let $v_i=\hat g_iv$.
Then $\hat g_i \hat p= \omega g_j\hat p$ where $\omega$ is a cubic
root of unity.
Let $g\in  \Gamma$ be a side pairing of the polyhedra obtained by
the development map.  We have
$
gg_ip=g_jp,
$
therefore $(g_j)^{-1}gg_ip=p$ and we conclude that $(g_j)^{-1}gg_i$ is parabolic
and then its lift $(\hat g_j)^{-1}\hat g\hat g_i$ is a parabolic
element (or the identity element) fixing $v$ up to a cubic root of unity.
%This implies that $\hat g\hat g_iv =\pm \hat g_j v$
We conclude that
$$
\hat g v_i= v_j.
$$
up to multiplication by a cubic root of unity.
We do that for each cycle of vertices and
obtain that the lift is special.
 \EPf

\begin{thm}\label{B(C)}
Let $M$ be a $3$-dimensional spherical CR manifold with a CR triangulation,
say $M=\Delta_1\cup \cdots \cup \Delta_k$ with each $\Delta_i$ a CR
tetrahedron. Suppose $M$ is non-compact with purely
parabolic boundary.
Then its invariant $\beta(M)\in B(\C)$, that is,
$$
\sum_i \delta(\beta(\Delta_i))=0.
$$
\end{thm}
\Pf The proof follows from the previous proposition.  We choose a
special lift and each face of a fixed tetrahedron appears four times
in the sum but the terms corresponding to a common face between two
tetrahedra appear with opposite sign. The only problem might arise
with the side pairing maps, but the surface terms obtained will
differ by terms (which are null) of the form $\omega\wedge a$ where
$\omega$ is a cubic root of unity. \EPf

The arguments  above are valid if we substitute $\C$ for $k=k_{_{\bf X}}$
where $k_{_{\bf X}}$ is the field generated by the cross-ratios of
a triangulation.  We will make the hypothesis that $k$ is a number field
in order to obtain a special lift:

\begin{lem}\label{special}
If a CR triangulation with $k=k_{_{\bf X}}$ a number field is such that
each vertex of the developed triangulation
is the fixed point of a   parabolic or elliptic element
(that is, we exclude loxodromic elements) of the holonomy group
then a finite cover of the triangulation has development with a special lift.
\end{lem}
Here the root of unity appearing in the definition of a special lift  is in $k$.
\vspace{.5cm}

\Pf
We first consider a normalized development of the vertices.
 By Proposition \ref{coordinates} the coordinates
of a normalised development of the triangulation are contained
in the invariant field.
  By Lemma 3.4 of \cite{Ge} any element $M$ in the holonomy
is such that $M^3\in SU(2,1,k)$. Let $\Gamma$ be the group generated by the
cubes.  It is a finite index subgroup of the holonomy group and we consider
the corresponding finite cover of the triangulation.  It defines the same
field $k$ but now the side pairings have lifts in $SU(2,1,k)$.
  Now we argue as before:
 we divide the vertices into classes, each one
obtained by translating a fixed vertex by elements of the holonomy
group.   Fixing a vertex $p\in S^3$ we consider all vertices
 identified to $p$ by the holonomy group.
That is, the
other vertices $p_i$
are obtained by $p_i=g_ip$ for $g_i\in \Gamma$. Suppose
$g_ip=g_jp$. Then $g_j^{-1}g_ip=p$ and by hypothesis $g_j^{-1}g_i$
is parabolic or elliptic with coefficients in $k$.
 Without loss of generality, we suppose that $p=\infty$ and
let
$\hat g_i\in SU(2,1,k)$ be a lift of $g_i$ and $v_i=\hat g_iv$.
Then $\hat g_i \hat p= \mu g_j\hat p$ where $\mu$ is a root of unity in the field $k$.
Let $g\in  \Gamma$ be a side pairing of the polyhedra obtained by
the development map.  We have
$
gg_ip=g_jp,
$
therefore $(g_j)^{-1}gg_ip=p$ and we conclude that $(g_j)^{-1}gg_i$ is parabolic
or elliptic
and then its lift $(\hat g_j)^{-1}\hat g\hat g_i$
 fixes $v$ up to a root of unity.
We conclude that
$$
\hat g v_i= v_j.
$$
up to multiplication by a root of unity.
We do that for each cycle of vertices and
obtain that the lift is special.
 \EPf

\begin{thm}\label{B(k)}
Let $k=k_{_{\bf X}}$ (which we suppose to be a number field)
 be the invariant field of a CR-triangulation of $M$
with  parabolic or elliptic boundary holonomy.
Then there exists an integer $d\geq 1$ such that
$d\beta(M)\in \mathcal{B}(k)$.
\end{thm}
\Pf
By the proposition above we obtain a special normalized lift of a certain
finite cover of $M$.  The proof then follows as in theorem \ref{B(C)} for the
cover so that $d\beta(M)\in \mathcal{B}(k)$ for an integer $d\geq 1$.
%We conclude by  using the fact that $\Lambda^2k^*$ is a free abelian group.
\EPf

\begin{thm}\label{imaginary}
Let $k_{_{\bf X}}$ be the  field of a CR-triangulation $(T,{\bf X})$.
Suppose $k_{\bf X}\subset k$ where $k$ is a purely imaginary quadratic
 extension of a totally real field.
If $\beta(T,{\bf X})\in {\mathcal{B}(k)}$ then it is torsion.
\end{thm}
\Pf We proved that $D(\beta(T,{\bf X}))=0$ and the result follows from
Borel's theorem.
\EPf

Let $M$ be a CR-triangulation.
Assume that its invariant $\beta(M)\in \mathcal{B}(\C)$.
Then we have $\rho(\beta(M))\in \C/\pi^2\Q$.
By Theorem \ref{vol}, we know $\Im \rho(\beta(M))=0$.

\begin{dfn}
Let $M$ be as above with $\beta(M)\in \mathcal{B}(\C)$. We define its Chern-Simons invariant to be the real part of
$\rho(\beta(M))$, denoted by $CS(M)$.
\end{dfn}

{\bf Remark} Theorem \ref{imaginary} implies that if the invariant field
associated to a CR structure is an imaginary quadratic extension of a totally
real field then $CS(M)=0$.

\section{Examples}
Besides the 5-term relation (\ref{5term}) at the beginning,
it is known (\cite{S1}) that we have two more identities in
$\mathcal{P}(\C)$:
\begin{equation}\label{5term1}
[z]+[z^{-1}]=0,
\end{equation}
\begin{equation}\label{5term2}
[z]+[1-z]=0.
\end{equation}
In the following, for a comlex number $a$, we will denote its complex conjugation by $\overline{a}$.

\subsection{Figure 8 Knot Complement}\label{figure8}
The Figure 8 knot complement $K$ can be glued by two ideal CR
tetrahedra.
Solving the equations in \cite{F2} imposing that $R(H_2)$ be real,
 we obtain the following one real parameter family of solutions:
Let
$$
\alpha = \sqrt{ 2- 4\beta^2+2\sqrt{5-8\beta}}.
$$
Observe that $\alpha$ is real for $I=\{ \beta \ |\ -\frac{1}{2}-
\frac{\sqrt{5}}{2}\leq \beta\leq \frac{5}{8}\ \}$.
Then one branch of solutions is given (for $\beta\in I$), by:
$$
  w_{12}=\beta+\frac{\alpha}{2}i;\;\; w_{21}=\beta-\frac{\alpha}{2}i;\;\;
  w_{34}=\beta+\frac{\alpha}{2}i;\;\; w_{43}=\beta-\frac{\alpha}{2}i.
$$
and
$$
  z_{12}=\frac{\sqrt{5-8\beta}-2\beta+1+{\alpha}i}{\sqrt{5-8\beta}+3-4\beta}
;\;\; z_{21}=\overline {z_{12}};\;\;
  z_{34}= z_{12};\;\; z_{43}=\overline {z_{12}}.
$$

\begin{lem}
The invariant
$$
[w_{12}]+[w_{21}]+[w_{34}]+[w_{43}]+[z_{12}]+[z_{21}]+[z_{34}]+[z_{43}]=0 \;\; \text{in} \;\; \mathcal {P}({\C}).
$$
\end{lem}

\Pf By the definition, we see that
LHS$=2([w_{12}]+[\overline{w_{12}}]+[z_{12}]+[\overline{z_{12}}])$.
We will check that
$$
  \frac{1}{1-\frac{1}{w_{12}}}=\frac{w_{12}}{w_{12}-1}=\overline {z_{12}}, \;\text{and} \;
  \frac{1}{1-\frac{1}{\overline{w_{12}}}}=\frac{\overline{w_{12}}}{\overline{w_{12}}-1}=z_{12}.
$$
Then it will follow that $[w_{12}]+[\overline{z_{12}}]=0$, and
$[\overline{w_{12}}]+[z_{12}]=0$. The lemma will be proved. We only
need to check the first one since the second is the complex
conjugation of the first one.

$$
\frac{w_{12}}{w_{12}-1}=\frac{\beta+\frac{\alpha}{2}i}{(\beta-1)+\frac{\alpha}{2}i}
=\frac{(\beta+\frac{\alpha}{2}i)((\beta-1)-\frac{\alpha}{2}i)}{(\beta-1)^2+\frac{\alpha^2}{4}}
=\frac{[\beta(\beta-1)+\frac{\alpha^2}{4}-\frac{\alpha}{2}i]}{(\beta-1)^2+\frac{\alpha^2}{4}}
$$
Multiplying by $2$ on the numerator and denominator of the last
fraction, we get
$$
\frac{w_{12}}{w_{12}-1}=\frac{2\beta(\beta-1)+\frac{\alpha^2}{2}- \alpha i}{2(\beta-1)^2+\frac{\alpha^2}{2}}.
$$
By the equation of $\alpha$ and $\beta$ at the beginning, we have
$$
2\beta^2+\frac{\alpha^2}{2}=1+\sqrt{5-8\beta}.
$$
Plug this in and simplify, we have
$$
\frac{w_{12}}{w_{12}-1}=\frac{(\sqrt{5-8\beta}-2\beta+1)-{\alpha}i}{\sqrt{5-8\beta}+3-4\beta}
=\overline{z_{12}}.
$$
\EPf

For $\beta=1/2$ we obtain $\alpha=\sqrt{3}$ and therefore the invariant
$$
  \beta(K)=4([\omega]+[\overline{\omega}]),
$$
where $\omega=\exp(\frac{2\pi i}{3})$ is a primitive cube root of
unity. Since $\overline{\omega}=\omega^{-1}$, by (\ref{5term1}), we
see in $\mathcal {P}({\C})$
$$
  [\omega]+[\overline{\omega}]=0.
$$
So $\beta(K)=0$.

In \cite[~Page 94]{F2}, there are two other representations from the
fundamental group of the figure eight complement to $PU(2,1)$. They have cyclic holonomy on the boundary. The
first invariants are
$$
  w_{12}=\frac{3}{8}+\frac{\sqrt{7}}{8}i;\;\; w_{21}=\frac{5}{4}+\frac{\sqrt{7}}{4}i;\;\;
  w_{34}=\frac{3}{8}-\frac{\sqrt{7}}{8}i;\;\; w_{43}=\frac{5}{4}-\frac{\sqrt{7}}{4}i.
$$
The second invariants are
$$
  t_{12}=\frac{3}{2}+\frac{\sqrt{7}}{2}i;\;\; t_{21}=-\frac{1}{4}-\frac{\sqrt{7}}{4}i;\;\;
  t_{34}=\frac{3}{2}-\frac{\sqrt{7}}{2}i;\;\; t_{43}=-\frac{1}{4}+\frac{\sqrt{7}}{4}i;\;\;
$$
Let $F=\Q(\sqrt{-7})$.  Then we have the invariants
\[
\beta_1(K)=2([w_{12}]+[w_{21}]+[w_{34}]+[w_{43}]), \;\; \beta_2(K)=2([t_{12}]+[t_{21}]+[t_{34}]+[t_{43}])\in \mathcal {P}(F)
\]

\begin{prop}
(1). $\beta_1(K)=-2c_{F}$ in $\mathcal {P}(F)$. \\
(2). $\beta_2(K)=0$ in $\mathcal {P}(F)$. \\
(3). $\beta_1(K)\in \mathcal{B}(F)$ is a non-trivial torsion of order $3$.
\end{prop}

\Pf
(1). Let's put
$$
 a=w_{34}=\frac{3}{8}-\frac{\sqrt{7}}{8}i, \;\; b=w_{43}=\frac{5}{4}-\frac{\sqrt{7}}{4}i,
$$
then
$$
 \frac{1}{2}b=\frac{5}{8}-\frac{\sqrt{7}}{8}i=1-\overline{a},\;\;\;
 \frac{1}{2}\overline{b}=\frac{5}{8}+\frac{\sqrt{7}}{8}i=1-a,
$$
and
$$
 [\overline{a}]=c_F-[1-\overline{a}]=c_F-[\frac{1}{2}b],\;\;\; [a]=c_F-[1-a]=c_F-[\frac{1}{2}\overline{b}].
$$
Therefore,
$$
  \beta_1(W)=2([b]+[\overline{b}]-[\frac{1}{2}b]-[\frac{1}{2}\overline{b}]+2c_F).
$$
By the five-term equation (\ref{5term}), take $x=\frac{1}{2}$,
 $y=\frac{1}{2}b$, then we obtain
$$
 [\frac{1}{2}]-[\frac{1}{2}b]+[b]-[\frac{1-2}{1-(\frac{1}{2}b)^{-1}}]+[\frac{1-\frac{1}{2}}{1-\frac{1}{2}b}]=0.
$$
Direct computations show that
$$
  s=\frac{1-\frac{1}{2}}{1-\frac{1}{2}b}=\frac{3}{4}-\frac{\sqrt{7}}{4}i,\;\; \text{and} \;\; \frac{1-2}{1-(\frac{1}{2}b)^{-1}}=\frac{1}{\frac{1}{4}+\frac{\sqrt{7}}{4}i}=\frac{1}{1-s}.
$$
Hence,
$$
[b]-[\frac{1}{2}b]= [\frac{1}{1-s}]-[s]-[\frac{1}{2}].
$$
Similarly by taking $x=\frac{1}{2}$, $y=\frac{1}{2}\overline{b}$, we have
$$
[\overline{b}]-[\frac{1}{2}\overline{b}]= [\frac{1}{1-\overline{s}}]-[\overline{s}]-[\frac{1}{2}].
$$
Hence,
$$
 \beta_1(K)=2([\frac{1}{1-s}]-[s]+[\frac{1}{1-\overline{s}}]-[\overline{s}]-c_F)+4c_F
          =2([\frac{1}{1-s}]-[s]+[\frac{1}{1-\overline{s}}]-[\overline{s}])+2c_F.
$$
By Lemma 1.2 of \cite{S1}, in $\mathcal {P}(F)$ we have
$$
  2([\frac{1}{1-s}])=-2[1-s],  \;\; \text{and} \;\;   2([\frac{1}{1-\overline{s}}])=-2[1-\overline{s}],
$$
so
$$
  \beta_1(K)=-2c_F-2c_F+2c_F=-2c_F.
$$

(2). Notice that
$$
  t_{43}=1-w_{43},\;\; t_{21}=1-w_{21},\;\; t_{12}=w_{34}^{-1};\;\; t_{34}=w_{12}^{-1}.
$$
Hence we have
$$
[t_{43}]=c_F-[w_{43}],\;\; [t_{21}]=c_F-[w_{21}], \;\; 2[t_{12}]=-2[w_{34}], \;\; 2[t_{34}]=-2[w_{12}],
$$
and
$$
 \beta_2(K)=2([t_{12}]+[t_{21}]+[t_{34}]+[t_{43}])=4c_F-\beta_1(K)=6c_F=0.
$$
The last equality comes from \cite[~Lemma 1.5(a)]{S1}.\\

(3). Since $c_F\in \mathcal{B}(F)$, by (1), $\beta_1(K)\in \mathcal{B}(F)$. It is easy to see that for $F=\Q(\sqrt{-7})$, $\mu(F)=\{\pm 1\}=\Z/2$. Hence $\text{Tor}(\mu(F),\mu(F))=\Z/2$. The result follows by Lemma \ref{blem2} .

\EPf
\begin{cor}
Both $\beta_1(K)$ and $\beta_2(K)$ are zero in $\mathcal{B}(\C)$.
\end{cor}

\Pf
$\beta_2(K)=0$ in $\mathcal{B}(\C)$ since it is already zero in $\mathcal {P}(F)$. For $\beta_1(K)$, since it is a torsion in $\mathcal{B}(F)$, it will be a torsion in $\mathcal{B}(\C)$. We know that $\mathcal{B}(\C)$
is torsion-free, hence $\beta_1(K)=0$ in $\mathcal{B}(\C)$.
\EPf
\subsection{Whitehead Link Complement}\label{whitehead}

The Whitehead link complement $W$ can be glued by four ideal CR
tetrahedra.  See \cite{Sc, Ge}. Their CR invariants are:
$$
   A_{01}=-\frac{1}{8}+\frac{\sqrt{15}}{8}i;\;\; A_{10}=-2; \;\;
   A_{23}=-\frac{3}{4}+\frac{\sqrt{15}}{4}i;\;\;
   A_{32}=\frac{1}{2}-\frac{\sqrt{15}}{6}i.
$$

$$
   B_{01}=-2;\;\; B_{10}=\overline{A}_{01};\;\;
   B_{23}=\overline{A}_{32};\;\; B_{32}=\overline{A}_{23}.
$$

$$
  C_{01}=\overline{A}_{01};\;\; C_{10}=-2;\;\;
  C_{23}=\overline{A}_{23};\;\; C_{32}=\overline{A}_{32}.
$$

$$
  D_{01}=-2;\;\; D_{10}=A_{01};\;\; D_{23}=A_{32};\;\;
  D_{32}=A_{23}.
$$

Let $F=\Q(\sqrt{-15})$.  Then we have
\[
\beta(W)=4[-2]+2([A_{01}]+[\overline{A}_{01}]+[A_{23}]+[\overline{A}_{23}]+[A_{32}]+[\overline{A}_{32}])\in
\mathcal {P}(F)
\]

In $\mathcal {P}(F)$, by \cite[~Lemma 1.2]{S1}, for any $z\ne 0,1$ we have :
\begin{equation}\label{2tor}
2([z]+[z^{-1}])=0.
\end{equation}

\begin{lem}\label{6lem}
Let $z\in F-\{0,1,-1\}$. Then we have
 $$
   2[z^2]=4[z]+4[-z].
 $$
\end{lem}
\Pf
By \cite[~Lemma 4.5]{S1}, we have $[z^2]=2([z]+[-z]+[-1])$. By (\ref{2tor}), we have $4[-1]=0$. Hence
$$
  2[z^2]=4[z]+4[-z]+4[-1]= 2[z^2]=4[z]+4[-z].
$$
\EPf

\begin{prop}
(1). $\beta(W)=4[\frac{1}{2}]$ in $\mathcal {P}(F)$. \\
(2). $\beta(W)\wedge (1-\beta(W))=0$, hence $\beta(W)\in \mathcal{B}(F)$.\\
(3). $\beta(W)$ has order $3$ in $\mathcal{B}(F)$.
\end{prop}

\Pf
(1). Since
 $1-A_{32}=\frac{1}{2}+\frac{\sqrt{15}}{6}i=\overline{A}_{32}$, by
 Definition \ref{bcf}, we get
 \[
   [A_{32}]+[\overline{A}_{32}]=[A_{32}]+[1-A_{32}]=c_F.
 \]
 In the 5-term equation (\ref{5term}), take $x=A_{01}$,
 $y=A_{01}\overline{A}_{01}=|A_{01}|^2$, we obtain in $\mathcal
 {P}(F)$:
 \[
   [A_{01}]-[|A_{01}|^2]+[\overline{A}_{01}]-[\frac{1-A_{01}^{-1}}{1-|A_{01}|^{-2}}]+[\frac{1-A_{01}}{1-|A_{01}|^{2}}]=0.
 \]
 Therefore,
 \[
  [A_{01}]+[\overline{A}_{01}]=[|A_{01}|^2]+[\frac{1-A_{01}^{-1}}{1-|A_{01}|^{-2}}]-[\frac{1-A_{01}}{1-|A_{01}|^{2}}].
 \]
 Direct computations show that
 \[
   \frac{1-A_{01}^{-1}}{1-|A_{01}|^{-2}}=A_{23}^{-1}.
 \]

 \[
   \frac{1-A_{01}}{1-|A_{01}|^{2}}=1- \frac{1}{\overline{A}_{23}}.
 \]

 Therefore,
 \[
   [A_{01}]+[\overline{A}_{01}]=[\frac{1}{4}]+[
   A_{23}^{-1}]-[1- \frac{1}{\overline{A}_{23}}].
 \]
We obtain that
\begin{alignat}{2}
 \beta(W)&=4[-2]+2([A_{01}]+[\overline{A}_{01}]+[A_{23}]+[\overline{A}_{23}]+[A_{32}]+[\overline{A}_{32}])\notag \\
        &=4[-2]+2([\frac{1}{4}]+[A_{23}^{-1}]-[1- \frac{1}{\overline{A}_{23}}]+[A_{23}]+[\overline{A}_{23}]+c_F)
        \notag \\
        &=4[-2]+2([\frac{1}{4}]+[A_{23}^{-1}]-[1-\frac{1}{\overline{A}_{23}}]-[\frac{1}{\overline{A}_{23}}]+[\frac{1}{\overline{A}_{23}}]+[A_{23}]+[\overline{A}_{23}]+c_F)
        \notag \\
        &=4[-2]+2([\frac{1}{4}]+[A_{23}^{-1}]-c_F+[\frac{1}{\overline{A}_{23}}]+[A_{23}]+[\overline{A}_{23}]+c_F)
        \notag \\
        &=4[-2]+2([\frac{1}{4}]+[A_{23}]+[A_{23}^{-1}]+[\overline{A}_{23}]+[(\overline{A}_{23})^{-1}]) \notag \\
        &=4[-2]+2[\frac{1}{4}] \;\;\; (\text{By}\; (\ref{2tor}))\notag
\end{alignat}

Next by Lemma \ref{6lem} and (\ref{2tor}), we have in $\mathcal {P}(F)$
\[
   2[z^2]=4[z]+4[-z]\;\; \text{and} \;\; 2[z]=-2[z^{-1}]
\]
Therefore,
\[
  2[\frac{1}{4}]=4[\frac{1}{2}]+4[-\frac{1}{2}]=4[\frac{1}{2}]-4[-2].
\]
Now in $\mathcal {P}(F)$ we have
\[
  \beta(W)=4[-2]+4[\frac{1}{2}]-4[-2]=4[\frac{1}{2}].
\]

(2). Since
$$
  \delta(\rho(W)=4(\frac{1}{2} \wedge (1-\frac{1}{2}))=4(\frac{1}{2} \wedge \frac{1}{2})=0,
$$
$\beta(W)\in \mathcal{B}(F)$.

(3). Notice $\beta(W)=4[\frac{1}{2}]=2([\frac{1}{2}]+[1-\frac{1}{2}])=2c_F$. It is easy to see that for
$F=\Q(\sqrt{-15})$, $\mu(F)=\{\pm 1\}=\Z/2$. Hence $\text{Tor}(\mu(F),\mu(F))=\Z/2$ and it has no element of order
$3$. By Lemma \ref{blem2}, $\beta(W)=2c_F$ is a nonzero $3$-torsion in $\mathcal{B}(F)$.
\EPf

\begin{cor}
$\beta(W)=0$ in $\mathcal{B}(\C)$.
\end{cor}

\Pf
It follows from the fact that $\mathcal{B}(\C)$ is torsion-free.
\EPf

\section*{Acknowledgements} E. Falbel thanks Francis Brown, Herbert Gangl, Julien Grivaux and Julien March\'{e} for several helpful conversions. Part of the work was done when Q. Wang visited University Paris VI in January 2010. He thanks University Paris VI for hospitality and support. Q. Wang is grateful to Stephen Lichtenbaum for stimulating discussions.

%%%%%%%%%%%%%%%%%%%%%%%%%%%%%%%

%%%%%%%%%%%%%%%%%%%%%%%%%%%%%%%


\begin{thebibliography}{ZZ99}

\bibitem[A]{A}
J. W. Alexander ; The combinatorial theory of complexes.  Ann. of Math. (2)  31  (1930),  no. 2, 292--320.

\bibitem[B]{B} S. Bloch ; Applications of the dilogarithm function in algebraic K-theory and algebraic geometry, in:
    Proc. of the International Symp. on Alg. Geom., Tokyo, 1978.

\bibitem[B1]{B1} S. Bloch ; Higher regulators, algebraic $K$-theory, and zeta functions of elliptic curves. CRM
    Monograph Series, 11. American Mathematical Society, Providence, RI, 2000.

\bibitem[Bo]{Bo} A. Borel ;
Cohomologie de ${\rm SL}_{n}$ et valeurs de fonctions zeta aux points entiers.
Ann. Scuola Norm. Sup. Pisa Cl. Sci. (4)  4  (1977), no. 4, 613--636.

%\bibitem[Bu]{Bu} J. I. Burgos Gil ; The regulators of Beilinson and Borel. CRM Monograph Series, 15. American Mathematical Society, Providence, RI, 2002.

\bibitem[BS]{BS}  D. Burns, S. Shnider ;  Spherical Hypersurfaces
in Complex Manifolds.  Invent. Math. 33 (1976), 223-246.


\bibitem[C]{C} E. Cartan ; { Sur le groupe de la g\'eom\'etrie
hypersph\'erique}, Comm. Math. Helv. 4 (1932), 158-171.

\bibitem[Ca]{Ca}
B. G. Casler ; An imbedding theorem for connected $3$-manifolds with boundary.  Proc. Amer. Math. Soc.  16  1965
559--566.

\bibitem[CS]{CS}
S. S. Chern, J. Simons ; Characteristic forms and geometric invariants. Ann. of Math. (2) 99 (1974), 48--69.

%\bibitem[Cu]{Cu}  M. Culler ; Lifting Representations to Covering Groups, Advances in Mathematics 59 (1986), 64-70.

%\bibitem[CS]{CS}  M. Culler, P. B. Shalen ; Varieties of group representations and splittings of 3-manifolds, Ann. of Math. 117 (1983) 109-146.

\bibitem[D]{D} J. Dupont ; {The dilogarithm as a characteristic
class for flat bundles}, Journal of pure and applied algebra 44
(1987) 137-164.

\bibitem[DS]{DS}  J. Dupont; C. H. Sah; Scissors congruences. II. J. Pure Appl. Algebra 25 (1982), no. 2, 159--195.

\bibitem[DZ]{DZ}  J. Dupont, C. K. Zickert ;
A dilogarithm formula for the Cheeger-Chern-Simons class.
Geometry and Topology 10 (2006) 1347-1372.

\bibitem[F1]{F1} E. Falbel ; {Geometric structures associated to triangulations as fixed point sets of involutions}.
Topology and its Applications  154  (2007), no. 6, 1041-1052.

\bibitem[F2]{F2} E. Falbel ; { A spherical CR structure on the
complement of the figure eight knot with discrete holonomy}.
Journal of Differential Geometry 79 (2008) 69-110.

\bibitem[F3]{F3} E. Falbel ;
A volume function for Spherical CR Tetrahedra.
To appear in Quarterly Journal of Mathematics.


\bibitem[FP]{FP} E. Falbel, I. D. Platis ; The $PU(2,1)$ configuration space of four points in $S\sp 3$ and the
    cross-ratio variety.  Math. Ann.  340  (2008),  no. 4, 935--962.

\bibitem[Ge1]{Ge1} J. Genzmer ; Trace fields of subgroups of $SU(n,1)$. Preprint
2009.

\bibitem[Ge]{Ge} J. Genzmer ; Sur les triangulations des structures CR
sph\'eriques. Thesis, Paris VI (2010).


%\bibitem[Gu]{Gu} A. Guichardet ; Cohomologie des groupes topologiques et des alg\`ebres de Lie. Textes Math\'ematiques 2. CEDIC, Paris, 1980.

\bibitem[G]{G} W. M. ~Goldman ; Complex Hyperbolic Geometry.
Oxford Mathematical Monographs. Oxford University Press (1999).


\bibitem[J]{J} H. ~Jacobowitz ; An Introduction to CR Structures.
Mathematical Surveys and Monographs {\bf 32}, American Math.\ Soc.\ (1990).

%\bibitem[KR]{KR} A. Kor\' anyi, H. M. Reimann ; The complex cross ratio on the Heisenberg group. Enseign. Math. (2) 33 (1987), no. 3-4, 291--300.

\bibitem[Ma]{Ma}
S. V. Matveev ; Transformations of special spines, and the Zeeman conjecture.  Izv. Akad. Nauk SSSR Ser. Mat.  51
(1987),  no. 5, 1104--1116, 1119;  translation in  Math. USSR-Izv.  31  (1988),  no. 2, 423--434.

\bibitem[Ma1]{Ma1}
S. V. Matveev ; Algorithmic topology and classification of 3-manifolds.
 Algorithms and Computation in Mathematics, 9. Springer, Berlin, 2007.

\bibitem[Mc]{Mc} McReynolds, Ben ; Arithmetic PU(2,1)...


\bibitem[M]{M} J. Milnor ; { Hyperbolic geometry: the first 150 years.},
Bull. Am.  Math. Soc. 6,  (1982), 9-24.

\bibitem[NR]{NR} W. Neumann, A. Reid ; Arithmetic of hyperbolic manifolds. Topology '90 (Columbus, OH, 1990),
    273--310, Ohio State Univ. Math. Res. Inst. Publ., 1, de Gruyter, Berlin, 1992.



\bibitem[NZ]{NZ} W. Neumann, D. Zagier ; Volumes of hyperbolic three-manifolds.  Topology  24  (1985),  no. 3,
    307--332.


\bibitem[NY1]{NY1} W. Neumann, Jun Yang ;
Rationality problems for $K$-theory and Chern-Simons invariants of hyperbolic $3$-manifolds.
Enseign. Math. (2) 41 (1995), no. 3-4, 281--296.

\bibitem[NY]{NY} W. Neumann, Jun Yang ;
 Bloch invariants of hyperbolic 3-manifolds.
Duke Math. Journal 96 (1999) 25-59.

%\bibitem[O]{O} J. Osterl\'e ; Polylogarithmes. S\'eminaire Bourbaki, exp. 762, 1992-1993.
\bibitem[P]{P}
R. Piergallini ; Standard moves for standard polyhedra and spines. Third National Conference on Topology (Italian)
(Trieste, 1986).  Rend. Circ. Mat. Palermo (2) Suppl.  No. 18  (1988), 391--414.

\bibitem[PP]{PP} J. R. Parker, I. D. Platis ;  Complex hyperbolic Fenchel-Nielsen coordinates.  Topology  47  (2008),
    no. 2, 101--135.

\bibitem[PP1]{PP1} J. R. Parker, I. D. Platis ; Global, geometrical coordinates on Falbel's cross-ratio variety
Canadian Mathematical Bulletin 52 (2009) 285-294.

\bibitem[R]{R} A. Reid ;
A note on trace-fields of Kleinian groups.  Bull. London Math. Soc.  22  (1990),  no. 4, 349--352.

\bibitem[Sah]{Sah} C. H. Sah, Homology of classical Lie groups made discrete. III. J. Pure Appl. Algebra 56 (1989),
    no. 3,
269--312.

\bibitem[Sc]{Sc} R. E. Schwartz ; Spherical CR Geometry and Dehn Surgery.  Annals of Math. Studies vol. 165 (2007).

\bibitem[S1]{S1}A. A. Suslin; $K_3$ of a field and the Bloch group.
Proc. of the Steklov Institute of Math. Issue 4 (1991), 217-238.

\bibitem[S2]{S2}A. A. Suslin; The group $K_3$ for a field. Math. USSR. Izvestiya, Vol 36(1991), No.3, 541-565.


\bibitem[T]{T} W. Thurston ; The geometry and topology of 3-manifolds.
Lecture notes 1979.

%\bibitem[Z]{Z} D. Zagier ; The dilogarithm function. Frontiers in number theory, physics and geoemtry. II, 3-65, Springer, Berlin, 2007.

%\bibitem[W]{W} P. Will ; The punctured torus and Lagrangian triangle groups in ${\rm PU}(2,1)$.  J. Reine Angew. Math.  602  (2007), 95--121.

\bibitem[W1]{W1} P. Will ;
Traces, Cross-ratios and 2-generator Subgroups of PU(2,1).
To appear in Can. J. Math.

\bibitem[W2]{W2} P. Will ;
Bending Fuchsian representations of fundamental groups of cusped surfaces in PU(2,1). Preprint 2008.

\bibitem[Z]{Z} D. Zagier ; The dilogarithm function. Frontiers in number theory, physics and geometry. II, 3-65,
    Springer, Berlin, 2007.

\end{thebibliography}
\end{document}